\documentclass{article}
\usepackage[utf8]{inputenc}
\usepackage{amsmath}
\usepackage{amsthm}
\usepackage{amssymb}
\usepackage{amsfonts}
\usepackage{mathtools}
\usepackage{stmaryrd}
\usepackage{tikz-cd}
\usepackage{enumerate}
\usepackage{hyperref} \hypersetup{colorlinks=true,citecolor=[rgb]{0.153,0.255,0.545},urlcolor=black,colorlinks=true,linkcolor=[rgb]{0.153,0.255,0.545}}

\usepackage[style=numeric,doi=false,url=false,isbn=false]{biblatex}
\renewbibmacro{in:}{}

\usepackage[top=1in, bottom=1.25in, left=1.25in, right=1.25in]{geometry}

\newcommand{\ZZ}{\mathbb{Z}}

\newcommand{\CC}{\mathbb{C}}
\newcommand{\FF}{\mathbb{F}}
\newcommand{\PP}{\mathbb{P}}
\renewcommand{\AA}{\mathbb{A}}
\newcommand{\GG}{\mathbb{G}}

\newcommand{\K}{\textup{K}_0}
\newcommand{\id}{\textup{id}}
\newcommand{\AGL}{\textup{AGL}}
\newcommand{\SL}{\textup{SL}}
\newcommand{\GL}{\textup{GL}}
\newcommand{\PGL}{\textup{PGL}}

\newcommand{\point}{\star}
\newcommand{\Var}{\textup{\bf Var}}
\newcommand{\Mod}{\textup{\bf Mod}}
\DeclareMathOperator{\Spec}{Spec}

\newcommand{\Hom}{\textup{Hom}}

\newcommand{\Bd}{\textup{\bf Bd}}
\newcommand{\Bdp}{\textup{\bf Bdp}}

\newcommand{\Grpd}{\textup{\bf Grpd}}

\newcommand{\Span}{\textup{Span}}

\DeclareMathOperator{\tr}{tr}

\DeclareMathAlphabet\mathbfcal{OMS}{cmsy}{b}{n}

\newcommand{\smallbdmultiplication}{{\renewcommand{\bdscale}{0.25}\bdmultiplication}}

\newcommand{\smatrix}[1]{\left(\begin{smallmatrix}#1\end{smallmatrix}\right)}

\newtheorem{theorem}{Theorem}[section]
\newtheorem{maintheorem}{Theorem}

\newtheorem{proposition}[theorem]{Proposition}
\newtheorem{lemma}[theorem]{Lemma}

\theoremstyle{definition}

\newtheorem{remark}[theorem]{Remark}

\newtheorem{definition}[theorem]{Definition}
\numberwithin{theorem}{section}
\numberwithin{equation}{section}

\newcommand{\bdscale}{0.5}

\def\bdgenus {
\begin{tikzpicture}[semithick, scale=\bdscale, baseline=-0.5ex]
\begin{scope}
    \draw (-1,0) ellipse (0.2cm and 0.4cm);
    \draw (-1,0.4) .. controls (-0.5,0.6) and (0.5,0.6) .. (1,0.4);
    \draw (-1,-0.4) .. controls (-0.5,-0.6) and (0.5,-0.6) .. (1,-0.4);
    \draw (-0.5,0.1) .. controls (-0.5,-0.125) and (0.5,-0.125) .. (0.5,0.1);
    \draw (-0.4,0.0) .. controls (-0.4,0.0625) and (0.4,0.0625) .. (0.4,0.0);
    \draw (1,0) ellipse (0.2cm and 0.4cm);
\end{scope}
\end{tikzpicture}
}

\def\bdmultiplication {
\begin{tikzpicture}[semithick, scale=\bdscale, baseline=-0.5ex]
\begin{scope}
    \draw (-1,0) ellipse (0.2cm and 0.4cm);
    \draw (1,0.5) ellipse (0.2cm and 0.4cm);
    \draw (1,-0.5) ellipse (0.2cm and 0.4cm);
    \draw (-1,0.4) .. controls (0,0.4) and (0,0.9) .. (1,0.9);
    \draw (-1,-0.4) .. controls (0,-0.4) and (0,-0.9) .. (1,-0.9);
    \draw (1,0.1) .. controls (0.2,0.1) and (0.2,-0.1) .. (1,-0.1);
\end{scope}
\end{tikzpicture}
}

\def\bdprojplane {
\begin{tikzpicture}[semithick, scale=0.5*\bdscale, baseline=-0.5ex]
\begin{scope}
    \filldraw[fill=black!20!white, draw=black] (-2,-1) rectangle (2,1);
    \filldraw[fill=white, draw=black] (-1,0) circle (0.5);
    \filldraw[fill=white, draw=black] (1,0) circle (0.5);
    \draw [-{to[width=2mm,length=3mm]}] (0, 1) to ++(0.1, 0);
    \draw [-{to[width=2mm,length=3mm]}] (0, -1) to ++(-0.1, 0);
    \draw [-{to[width=2mm,length=3mm]}] (2, 0) to ++(0, -0.1);
    \draw [-{to[width=2mm,length=3mm]}] (-2, 0) to ++(0, 0.1);
\end{scope}
\end{tikzpicture}
}

\def\bdpunit {
\begin{tikzpicture}[semithick, scale=\bdscale, baseline=-0.5ex]
\begin{scope}
    \draw (0,0) ellipse (0.2cm and 0.4cm);
    \draw (0,-0.4) arc (-90:90:0.75cm and 0.4cm);
    \draw[black,fill=black] (0.2,0) circle (.4ex);
\end{scope}
\end{tikzpicture}
}

\def\bdpcounit {
\begin{tikzpicture}[semithick, scale=\bdscale, baseline=-0.5ex]
\begin{scope}
    \draw (0,0.4) arc (90:-90:0.2cm and 0.4cm);
    \draw (0,0.4) arc (90:270:0.75cm and 0.4cm);
    \draw[dashed] (0,0.4) arc (90:270:0.2cm and 0.4cm);
    \draw[black,fill=black] (0.2,0) circle (.4ex);
\end{scope}
\end{tikzpicture}
}

\def\bdpgenus {
\begin{tikzpicture}[semithick, scale=\bdscale, baseline=-0.5ex]
\begin{scope}
    \draw (-1,0) ellipse (0.2cm and 0.4cm);
    \draw (-1,0.4) .. controls (-0.5,0.45) and (0.5,0.45) .. (1,0.4);
    \draw (-1,-0.4) .. controls (-0.5,-0.45) and (0.5,-0.45) .. (1,-0.4);
    \draw (-0.5,0.1) .. controls (-0.5,-0.125) and (0.5,-0.125) .. (0.5,0.1);
    \draw (-0.4,0.0) .. controls (-0.4,0.0625) and (0.4,0.0625) .. (0.4,0.0);
    \draw (1,0.4) arc (90:-90:0.2cm and 0.4cm);
    \draw[dashed] (1,0.4) arc (90:270:0.2cm and 0.4cm);
    \draw[black,fill=black] (-0.8,0) circle (.4ex);
    \draw[black,fill=black] (1.2,0) circle (.4ex);
\end{scope}
\end{tikzpicture}
}

\def\bdpprojplane {
\begin{tikzpicture}[semithick, scale=0.5*\bdscale, baseline=-0.5ex]
\begin{scope}
    \filldraw[fill=black!20!white, draw=black] (-2,-1) rectangle (2,1);
    \filldraw[fill=white, draw=black] (-1,0) circle (0.5);
    \filldraw[fill=white, draw=black] (1,0) circle (0.5);
    \draw [-{to[width=3mm,length=3mm]}] (0, 1) to ++(0.1, 0);
    \draw [-{to[width=3mm,length=3mm]}] (0, -1) to ++(-0.1, 0);
    \draw [-{to[width=3mm,length=3mm]}] (2, 0) to ++(0, -0.1);
    \draw [-{to[width=3mm,length=3mm]}] (-2, 0) to ++(0, 0.1);
    \fill (-1, 0.5) circle (.8ex);
    \fill (1, 0.5) circle (.8ex);
\end{scope}
\end{tikzpicture}
}

\def\bdgenerictube#1 {
\begin{tikzpicture}[semithick, scale=\bdscale, baseline=-0.5ex]
\begin{scope}
    \draw (-1,0) ellipse (0.2cm and 0.4cm);
    \draw (-1,0.4) cos (-0.875, 0.5) sin (-0.75,0.6) cos (-0.625,0.5) sin (-0.5,0.4) cos (-0.375,0.5) sin (-0.25, 0.6) cos (-0.125,0.5) sin (0, 0.4) cos (0.125,0.5) sin (0.25,0.6) cos (0.375,0.5) sin (0.5,0.4) cos (0.625,0.5) sin (0.75,0.6) cos (0.875,0.5) sin (1,0.4);
    \draw (-1,-0.4) cos (-0.875, -0.5) sin (-0.75,-0.6) cos (-0.625,-0.5) sin (-0.5,-0.4) cos (-0.375,-0.5) sin (-0.25, -0.6) cos (-0.125,-0.5) sin (0,-0.4) cos (0.125,-0.5) sin (0.25,-0.6) cos (0.375,-0.5) sin (0.5,-0.4) cos (0.625,-0.5) sin (0.75,-0.6) cos (0.875,-0.5) sin (1,-0.4);
    \draw (1,0) ellipse (0.2cm and 0.4cm);
    \draw (0,0) node {$#1$};
\end{scope}
\end{tikzpicture}
}

\title{\Large \textbf{A Topological Quantum Field Theory for \\ Character Varieties of Non-orientable Surfaces}}
\author{Jesse Vogel \\ \small Mathematical Institute, Leiden University, \href{mailto:j.t.vogel@math.leidenuniv.nl}{j.t.vogel@math.leidenuniv.nl}}
\date{}

\begin{document}

\maketitle

\vspace{-1cm}
\begin{abstract}
    In this paper, we study the $G$-representation and character varieties of non-orientable closed surfaces.
    By means of a geometric method based on a Topological Quantum Field Theory (TQFT), we compute the virtual classes of these varieties in the Grothendieck ring of varieties for $G$ equal to $\AGL_1$ and $\SL_2$. This method was already known and used in the case of orientable closed surfaces, and we extend it to the case of non-orientable surfaces.
    Furthermore, we provide a practical approach for explicitly computing the TQFT, allowing for more simplified and structured computations.
    Finally, for $G = \SL_2$ we describe and explain the relationship between the representation varieties of the orientable and non-orientable closed surfaces.
\end{abstract}

\section{Introduction}

Let $M$ be a closed connected manifold with finitely generated fundamental group $\pi_1(M)$, and $G$ an algebraic group over a field $k$. The set of group representations
\[ R_G(M) = \Hom(\pi_1(M), G) , \]
naturally carries the structure of an algebraic variety over $k$, and is called the \emph{$G$-representation variety} of $M$. Indeed, given a set of generators $\gamma_1, \ldots, \gamma_n$ of $\pi_1(M)$, the morphism
\[ R_G(M) \to G^n, \quad \rho \mapsto \left( \rho(\gamma_1), \ldots, \rho(\gamma_n) \right) \]
identifies the $G$-representation variety with a closed subvariety of $G^n$, and this structure can be shown to be independent of the chosen generators. The group $G$ acts on $R_G(M)$ by conjugation, identifying isomorphic representations, and taking the GIT quotient gives the \emph{$G$-character variety} of $M$,
\[ X_G(M) = R_G(M) \sslash G . \]
When $M = \Sigma_g$ is a closed surface of genus $g$, e.g.\ the underlying topological space of a smooth complex projective curve $C$, the corresponding $G$-representation variety can be expressed as
\begin{equation}
    \label{eq:representation_variety_orientable_surface}
    R_G(\Sigma_g) = \left\{ (A_1, B_1, \ldots, A_g, B_g) \in G^{2g} \;\bigg\vert\; \prod_{i = 1}^{g} [A_i, B_i] = 1 \right\} .
\end{equation}
The corresponding character variety plays an important role in non-abelian Hodge theory, being isomorphic to the moduli space of $G$-flat connections on $C$ and to the moduli space of semistable $G$-Higgs bundles on $C$, as real analytic spaces \cite{Simpson1992,Simpson1994a,Simpson1994b}. However, the complex structures on these three moduli spaces are different, making it very interesting to study the algebraic structures of these spaces, such as their mixed Hodge structure, or their virtual class in the Grothendieck ring of varieties.

In \cite{GonzalezLogaresMunoz2020}, Gonz\'alez-Prieto, Logares and Mu\~noz developed a geometric method to compute the virtual classes of the $G$-representation varieties $R_G(\Sigma_g)$ in the Grothendieck ring of varieties $\K(\Var_k)$ using a Topological Quantum Field Theory (TQFT). By a TQFT we mean a (lax) monoidal functor $Z : \Bd_n \to R\text{-}\Mod$ from the category of $n$-bordisms to the category of $R$-modules for some commutative ring $R$. One of the useful properties of a TQFT is that it assigns invariants to closed manifolds. Indeed, any closed manifold $M$ can be considered as a bordism $\varnothing \to \varnothing$, inducing a map $Z(M) : R \to R$ as $Z(\varnothing) = R$ by monoidality. This map is completely determined by the element $Z(M)(1) \in R$, which we refer to as the invariant associated to $M$. The TQFT that was constructed in \cite{GonzalezLogaresMunoz2020} associates to $M$ the virtual class $[ R_G(M) ]$ in $R = \K(\Var_k)$.

This method was used in \cite{Gonzalez2020} to compute the virtual classes of the parabolic $\SL_2$-representation and character varieties of $\Sigma_g$ over $k = \CC$ for any genus $g$. It was also used in \cite{HablicsekVogel2020,Vogel2023} to compute the virtual classes of $R_G(\Sigma_g)$ and $X_G(\Sigma_g)$ for $G$ the groups of upper triangular matrices of rank $2$, $3$, $4$ and $5$.
In this paper, we extend these results to the case of non-orientable surfaces with $G$ equal to $\AGL_1$ and $\SL_2$, the groups of affine transformations of the line and $2 \times 2$ matrices with determinant $1$, respectively.

By the classification of surfaces, any closed connected surface $\Sigma$ is either a sphere or a connected sum of tori or of projective planes. Let $N_r$ be the surface obtained by a connected sum of $r$ projective planes (that is, having demigenus $r$). Its fundamental group is given by
\[ \pi_1(N_r) = \Big\langle a_1, \ldots, a_r \;\Big\vert\; a_1^2 \cdots a_r^2 = 1 \Big\rangle , \]
and hence the corresponding $G$-representation variety can be expressed as
\begin{equation}
    \label{eq:representation_variety_nonorientable_surface}
    R_G(N_r) = \Big\{ (A_1, \ldots, A_r) \in G^{r} \;\Big\vert\; A_1^2 \cdots A_r^2 = 1 \Big\} .
\end{equation}
These representation varieties were also studied in \cite{LetellierRodriguezVillegas2022} by means of the so-called \emph{point-counting method}. This method, initiated by Hausel and Rodriguez-Villegas \cite{HauselRodriguezVillegas2008}, uses the character tables of $G(\FF_q)$ for finite fields $\FF_q$ to count the points of $R_G(\Sigma_g)$ (or $R_G(N_r)$) over $\FF_q$. If this count is polynomial in the order $q$ of the field, then it follows by Katz' theorem \cite[Theorem 6.1.2]{HauselRodriguezVillegas2008} that this polynomial is equal to the $E$-polynomial (also known as Serre polynomial) of the variety. Interestingly, for the surfaces $\Sigma_g$ only the dimension of the irreducible representations of $G(\FF_q)$ is needed, whereas for the surfaces $N_r$ also the Schur indicator of the irreducible representations is required.

Both the TQFT method and the point-counting method have advantages over each other. The invariant produced by the TQFT is more refined, that is, the $E$-polynomial of a complex variety is determined from its virtual class in the Grothendieck ring of varieties, but not vice versa. On the other hand, the point-counting method requires little computation once the character table of $G$ over finite fields is known, whereas the computations needed for the TQFT are often quite involved.


In Section \ref{sec:construction_TQFT} and \ref{sec:nonorientable_surfaces}, we recall the construction of the TQFT \cite{GonzalezLogaresMunoz2020} and extend it to non-orientable manifolds.
In Section \ref{sec:practical} we slightly modify the maps produced by the TQFT, resulting in a more practical approach to the explicit computations.
Finally, in Section \ref{sec:computations_AGL1} and \ref{sec:computations_SL2} we apply the method, explicitly computing the virtual classes of both the $G$-representation and character varieties of $N_r$ for $G$ equal to $\AGL_1$ and $\SL_2$.
In the case of $G = \SL_2$, we show that the orientable and non-orientable case are related through
\[ [R_{\SL_2}(\Sigma_g)] = [R_{\SL_2}(N_{2g})] \quad \textup{ for all } g \ge 0 , \]
and explain this relation from the point of view of the point-counting method.

Denote by $q$ the virtual class of the affine line $[\AA^1_k] \in \K(\Var_k)$. The main results are given by the following theorems.
\begin{maintheorem}
    \label{thm:virtual_class_nonorientable_AGL1}
    Let $k$ be any field of characteristic not equal to $2$. For any $r > 0$,
    \begin{enumerate}[(i)]
        \item the virtual class of the $\AGL_1$-representation variety of $N_r$ is
        \[ [R_{\AGL_1}(N_r)] = 2 q^{r - 1} \left(q - 1\right)^{r - 1} + q^{r - 1} \left(q - 1\right) , \]
        \item the virtual class of the $\AGL_1$-character variety of $N_r$ is
        \[ [X_{\AGL_1}(N_r)] = 2 (q - 1)^{r - 1} . \]
    \end{enumerate}
\end{maintheorem}

\begin{maintheorem}
    \label{thm:virtual_class_nonorientable_SL2}
    Let $k = \CC$. For any $r > 0$,
    \begin{enumerate}[(i)]
        \item the virtual class of the $\SL_2$-representation variety of $N_r$ is
        \begin{align*}
            [R_{\SL_2}(N_r)]
                &= \frac{1}{4} q^{r - 1} (q - 1) (q + 1)^{r - 1} ( (q - 1) (1 + (-1)^r) - (-2)^{r + 1} \\
                &\quad + \frac{1}{4} q^{r - 1} (q - 1)^{r - 1} (q + 1) ( (q - 1) (1 + (-1)^r) + 2^{r + 1} - 4 ) \\
                &\quad + q (q - 1)^{r - 1} (q + 1)^{r - 1}
                + q^{r - 1} (q - 1)^{r - 1} (q + 1)^{r - 1} ,
        \end{align*}
        \item the virtual class of the $\SL_2$-character variety of $N_r$ is
        \begin{align*}
            [X_{\SL_2}(N_r)]
                &= \frac{1}{4} q^{r - 2} (q - 1)^{r - 2} (((-1)^r + 1) (q - 1) + 2^{r + 1} - 12) \\
                &\quad + \frac{1}{4} q^{r - 2} (q + 1)^{r - 2} (((-1)^r + 1) (q - 1) - (-2)^{r + 1}) \\
                &\quad + q^{r - 2} (q - 1)^{r - 2} (q + 1)^{r - 2}
                + (q - 1)^{r - 2} (q + 1)^{r - 2}
                + q (q - 1)^{r - 2} \\
                &\quad + q (q + 1)^{r - 2}
                -2^{r - 1} q^{r - 2} ((-1)^r + 1) ,
        \end{align*}
    \end{enumerate}
    in the Grothendieck ring $\K^{\PP^1}(\Var_k)$.
\end{maintheorem}

The ring $\K^{\PP^1}(\Var_k)$ in the above theorem is the quotient of $\K(\Var_k)$ by the relations $[P] = [\PP^1_k] [X]$ for all $\PP^1$-fibrations $P \to X$ over $k$. Note that most useful invariants that factor through $\K(\Var_k)$, such as the $E$-polynomial, will also factor through $\K^{\PP^1}(\Var_k)$. We pass to this ring in order for the computations not to be overly complicated, but we strongly suspect the equalities in the theorem to still hold in $\K(\Var_k)$. 

\section{Construction of the TQFT}
\label{sec:construction_TQFT}

In this section, we briefly discuss the construction of the TQFT mentioned in the introduction. First of all, let us recall the definition of Grothendieck ring of varieties.

\begin{definition}
    Let $S$ be a variety over a field $k$. The \emph{Grothendieck ring of varieties} over $S$, denoted $\K(\Var/S)$, is the quotient of the free abelian group on the isomorphism classes of varieties over $S$, by relations of the form
    \[ [X] = [Z] + [X \setminus Z] \]
    for $Z \subset X$ a closed subvariety of a variety $X$ over $S$. Indeed $\K(\Var/S)$ carries a ring structure, where multiplication is given by
    \[ [X] \cdot [Y] = [X \times_S Y] \]
    for varieties $X$ and $Y$ over $S$, and the unit is $[S]$. The \emph{virtual class} of a variety $X$ over $S$ is its image in $\K(\Var/S)$. When $S = \Spec k$, we also denote the Grothendieck ring by $\K(\Var_k)$, and the virtual class of the affine line by $q = [\AA^1_k] \in \K(\Var_k)$.
\end{definition}

Let $G$ be an algebraic group over a field $k$. The TQFT that is used to compute the virtual classes of $G$-representation varieties in $\K(\Var_k)$ will be a lax monoidal functor $Z : \Bdp_n \to \K(\Var_k)\text{-}\Mod$ from the category of pointed $n$-bordisms to the category of $\K(\Var_k)$-modules. For a more elaborate construction of this TQFT, see \cite{Gonzalez2018,GonzalezLogaresMunoz2020,HablicsekVogel2020}. In these referenced works, all bordisms are assumed to be orientable. However, since our focus will be on non-orientable surfaces, we allow bordisms to be unoriented.

\begin{definition}
    Given two $(n - 1)$-dimensional closed manifolds $M_1$ and $M_2$, a \emph{bordism} from $M_1$ to $M_2$ is a compact $n$-dimensional manifold $W$ with boundary $\partial W$ together with inclusions
    \[ M_2 \overset{i_2}{\longrightarrow} W \overset{i_1}{\longleftarrow} M_1 \]
    such that $\partial W = i_1(M_1) \sqcup i_2(M_2)$. 
\end{definition}

\begin{definition}
    The \emph{category of pointed $n$-bordisms}, denoted $\Bdp_n$, is the category whose
    \begin{itemize}
        \item objects are pairs $(M, A)$, with $M$ a closed $(n - 1)$-dimensional manifold and $A \subset M$ a finite set of basepoints intersecting each connected component of $M$,
        \item morphisms $(M_1, A_2) \to (M_2, A_2)$ are equivalence classes of pairs $(W, A)$, with $W : M_1 \to M_2$ a bordism and $A \subset W$ a finite set of basepoints intersecting each connected component of $W$ such that $A \cap M_1 = A_1$ and $A \cap M_2 = A_2$. Two such pairs $(W, A)$ and $(W', A')$ are said to be equivalent if there is a diffeomorphism $F : W \to W'$ such that $F(A) = A'$ and such that the following diagram commutes:
        \begin{equation}
            \label{eq:diagram_def_bdp}
            \begin{tikzcd}[row sep=0em] & W \arrow{dd}{\wr} & \\ M_2 \arrow{ur} \arrow{dr} & & \arrow{ul} \arrow{dl} M_1 \\ & W' & \end{tikzcd}
        \end{equation}
        The composition of 1-morphisms $(W, A) : (M_1, A_1) \to (M_2, A_2)$ and $(W', A') : (M_2, A_2) \to (M_3, A_3)$ is $(W \sqcup_{M_2} W', A \cup A') : (M_1, A_1) \to (M_3, A_3)$. Although this operation is not well-defined on bordisms (there need not be a unique manifold structure on $W \sqcup_{M_2} W'$), it does define a well-defined operation on equivalence classes of bordisms \cite{Milnor1965}. Also note that unless $M = A = \varnothing$, there is not yet an identity morphism for $(M, A)$. For this reason we allow $(M, A)$ to be seen as a bordism from and to itself.
    \end{itemize}
    The category $\Bdp_n$ is a monoidal category, whose tensor product is given by taking disjoint unions, and whose unital object is the empty manifold $\varnothing$.
\end{definition}

Allowing for multiple basepoints directly leads to a generalized notion of the $G$-representation variety.
\begin{definition}
    Let $M$ be a compact manifold, possibly with boundary, and let $A \subset M$ be a finite set of basepoints. The \emph{$G$-representation variety} of $(M, A)$ is the variety
    \[ R_G(M, A) = \Hom_{\Grpd}(\Pi(M, A), G) , \]
    where $\Pi(M, A)$ denotes the fundamental groupoid of $(M, A)$.
\end{definition}

Let us elaborate a bit on the structure of $R_G(M, A)$. The groupoid $\mathcal{G} = \Pi(M, A)$ has finitely many connected components, where we say objects $a, b \in \mathcal{G}$ are \textit{connected} if $\Hom_{\mathcal{G}}(a, b)$ is non-empty.
Pick a subset $S = \{ a_1, \ldots, a_s \} \subset A$ such that each connected component of $\mathcal{G}$ contains exactly one of the $a_i$, and pick an arrow $f_a : a_i \to a$, with $a_i \in S$ in the same connected component as $a$, for each $a \in A \setminus S$.
Now a morphism of groupoids $\rho : \mathcal{G} \to G$ is uniquely determined by group morphisms $\rho_i : \mathcal{G}_i \to G$, where $\mathcal{G}_i = \Hom_{\mathcal{G}}(a_i, a_i)$, and a choice of $\rho(f_a) \in G$ for each $a \in A \setminus S$. Namely, any $\gamma : a \to b$ in $\mathcal{G}$ can be written as $\gamma = f_b \circ \gamma' \circ (f_a)^{-1}$ for some $\gamma' \in \mathcal{G}_i$, with $a_i \in S$ in the same connected component as $a$ and $b$, taking $f_{a_i} = \id$ if either $a$ or $b$ equals $a_i$. There is no restriction on the choice of $\rho(f_a)$, so we obtain
\begin{equation}
    \label{eq:the_hom_equation}
    R_G(M, A) = \Hom_{\Grpd}(\mathcal{G}, G) \cong \Hom(\mathcal{G}_1, G) \times \cdots \times \Hom(\mathcal{G}_s, G) \times G^{|A| - s} .
\end{equation}
Each of these factors naturally carries the structure of an algebraic variety. Namely, since $M$ is compact, it has the homotopy type of a finite CW-complex, so each $\mathcal{G}_i$ is finitely generated, and thus $\Hom(\mathcal{G}_i, G)$ can be identified with a subvariety of $G^m$ for some $m > 0$ as explained in the introduction. This gives $R_G(M, A)$ the structure of an algebraic variety, and this structure can be shown not to depend on the choices made. In particular, when $M$ is connected, we obtain
\begin{equation}
    \label{eq:G_repr_var_when_X_connected}
    R_G(M, A) \cong R_G(M) \times G^{|A| - 1} .
\end{equation}


Now, the TQFT $Z : \Bdp_n \to \K(\Var_k)\textup{-}\textbf{Mod}$ will be defined as the composition of two functors
\[ \Bdp_n \xrightarrow{\mathcal{F}} \Span(\Var_k) \xrightarrow{\mathcal{Q}} \K(\Var_k)\textup{-}\textbf{Mod} , \]
where $\Span(\Var_k)$ denotes the category of spans of varieties over $k$. Firstly, the functor $\mathcal{F}$, also known as the \emph{field theory}, is defined by
\begin{align*}
    \mathcal{F}(M, A) &= R_G(M, A) \\
    \mathcal{F}(W, A) &= \Big( R_G(M_1, A_1) \overset{i_1^*}{\longleftarrow} R_G(W, A) \overset{i_2^*}{\longrightarrow} R_G(M_2, A_2) \Big)
\end{align*}
for any object $(M, A)$ and bordism $(W, A) : (M_1, A_1) \to (M_2, A_2)$. Using the Seifert--van Kampen theorem for fundamental groupoids \cite{Brown1967}, one can show that $\mathcal{F}$ defines a functor, and is in fact a monoidal functor \cite{GonzalezLogaresMunoz2020}.

Secondly, the functor $\mathcal{Q}$, also known as the \emph{quantization functor}, is defined by
\begin{align*}
    \mathcal{Q}(X) &= \K(\Var/X) \\
    \mathcal{Q}(X \overset{f}{\longleftarrow} Z \overset{g}{\longrightarrow} Y) &= g_! \circ f^*
\end{align*}
for any variety $X$ over $k$ and span $X \overset{f}{\longleftarrow} Z \overset{g}{\longrightarrow} Y$ of varieties over $k$. Note that $\K(\Var/X)$ is indeed a $\K(\Var_k)$-module with
\[ [T] \cdot [W \xrightarrow{h} X] = [T \times_k W \xrightarrow{h \circ \pi_W} X] \]
for all varieties $T$ over $k$ and varieties $W$ over $X$, and that
\begin{align*}
    f^* &: \K(\Var/X) \to \K(\Var/Z), \quad [W \to X] \mapsto [W \times_X Z \to Z] \\
    g_! &: \K(\Var/Z) \to \K(\Var/Y), \quad [W \xrightarrow{h} Z] \mapsto [W \xrightarrow{g \circ h} Y]
\end{align*}
are indeed morphisms of $\K(\Var_k)$-modules. While $\mathcal{Q}$ is easily seen to be a functor, it is only lax monoidal, since the natural morphism
\[ \K(\Var/X) \otimes_{\K(\Var_k)} \K(\Var/Y) \to \K(\Var/(X \times_k Y)) \]
is usually not an isomorphism.

Finally, the composition of $\mathcal{F}$ and $\mathcal{Q}$ gives us the TQFT
\[ Z = \mathcal{Q} \circ \mathcal{F}: \Bdp_n \to \K(\Var_k)\text{-}\Mod , \]
which is a lax monoidal functor.
The crucial property of this TQFT is the following.
Let $M$ be a closed connected $n$-dimensional manifold, choose a basepoint $\point$ on $M$, and view $(M, \point)$ as a bordism $\varnothing \to \varnothing$. Then $\mathcal{F}(M, \point)$ is the span
\[ \begin{tikzcd} \point & \arrow[swap]{l}{t} R_G(M, \point) \arrow{r}{t} & \point \end{tikzcd} \]
and hence $Z(M, \point)(1) = t_! t^* (1) = [ R_G(M) ]$. This shows that the invariants associated to closed connected manifolds by this TQFT are precisely the classes of their $G$-representation varieties.

\section{Non-orientable surfaces}
\label{sec:nonorientable_surfaces}

We restrict our attention to the case of dimension $n = 2$. By the classification of surfaces, any closed orientable surface can be constructed as a composition of the following bordisms:
\begin{equation*}
    \renewcommand{\bdscale}{1.0}
    \bdpcounit : (S^1, \point) \to \varnothing \qquad \bdpgenus : (S^1, \point) \to (S^1, \point) \qquad \bdpunit : \varnothing \to (S^1, \point)
\end{equation*}
Indeed, a surface of genus $g$ can be written as $\bdpcounit \circ \bdpgenus^g \circ \bdpunit$ after adding sufficiently many basepoints.
Similarly, any closed non-orientable surface can be written as the connected sum of projective planes. For this reason, we consider the bordism
\begin{equation*}
    \renewcommand{\bdscale}{1.0}
    \bdpprojplane : (S^1, \point) \to (S^1, \point)
\end{equation*}
which is the real projective plane, with two holes and two basepoints, $a$ and $b$, seen as a bordism from $(S^1, a)$ to $(S^1, b)$.

The fundamental group of $S^1$ is $\pi_1(S^1) \cong \ZZ$, which implies $R_G(S^1, \point) \cong G$, and the fundamental group of the real projective plane with two punctures is the free group on two elements $\pi_1 \left( \bdprojplane, a \right) \cong F_2$. We choose generators $\alpha$ and $\beta$ of $\pi_1 \left( \bdpprojplane, a \right)$ as depicted in the following image, and a path $\gamma$ connecting $a$ and $b$.
\[ \begin{tikzpicture}[semithick, scale=3.0*\bdscale, baseline=-0.5ex]
\begin{scope}
    \filldraw[fill=black!20!white, draw=black] (-2,-1) rectangle (2,1);
    \filldraw[fill=white, draw=black] (-1,0) circle (0.5);
    \filldraw[fill=white, draw=black] (1,0) circle (0.5);
    \draw [-{To[width=3mm, length=3mm]}] (0, 1) to ++(0.1, 0);
    \draw [-{To[width=3mm, length=3mm]}] (0, -1) to ++(-0.1, 0);
    \draw [-{To[width=3mm, length=3mm]}] (2, 0) to ++(0, -0.1);
    \draw [-{To[width=3mm, length=3mm]}] (-2, 0) to ++(0, 0.1);
    \fill (-1, 0.5) circle (0.06);
    \fill (1, 0.5) circle (0.06);
    \draw [-{Latex[width=2mm, length=2mm]}] (-1, -0.5) to ++(0.1,0);
    \draw[dashed] (-2, -0.75) .. controls (-1.75, -0.75) and (-1.75, 0.5) .. (-1.5, 0.5) -- (-1, 0.5) .. controls (0, 0.5) and (0, 0.75) .. (2, 0.75);
    \draw [-{Latex[width=2mm, length=2mm]}] (1.6, 0.75) to ++(0.1,0);
    \draw (-1, 0.5) .. controls (0, 0.4) .. (1, 0.5);
    \draw [-{Latex[width=2mm, length=2mm]}] (0, 0.425) to ++(0.1, 0);
    \node at (-1, 0.3) {$a$};
    \node at (1, 0.3) {$b$};
    \node at (-1, -0.7) {$\beta$};
    \node at (1.6, 0.5) {$\alpha$};
    \node at (0, 0.2) {$\gamma$};
\end{scope}
\end{tikzpicture} \]

According to \eqref{eq:the_hom_equation}, we obtain $R_G \left( \bdprojplane, \{ a, b \} \right) \cong \Hom(F_2, G) \times G \cong G^3$. A generator for $\pi_1(S^1, b)$ is given by $\gamma \beta \alpha^2 \gamma^{-1}$, and hence the field theory for this bordism is given by the span
\[ \mathcal{F}\left(\bdpprojplane\right) = \left(\begin{tikzcd}[row sep=0em]
     G & \arrow[swap]{l}{u} G^3 \arrow{r}{v} & G \\
    B & \arrow[mapsto]{l} (B, A, C) \arrow[mapsto]{r} & C B A^2 C^{-1}
\end{tikzcd}\right) . \]
Consequently, the TQFT for this bordism is given by
 \[ Z \left( \bdpprojplane \right) = v_! \circ u^* : \K(\Var/G) \to \K(\Var/G) , \quad \left[ \begin{tikzcd} X \arrow{d}{f} \\ G \end{tikzcd} \right] \mapsto \left[ \begin{tikzcd}[column sep=0em] X \times G^2 \arrow{d} & (x, A, C) \arrow[mapsto]{d} \\ G & C f(x) A^2 C^{-1} \end{tikzcd} \right] . \]
Since the fundamental group of a disk is trivial, it is easy to see that the field theories for $\bdpunit$ and $\bdpcounit$ are given by
\[ \mathcal{F}\left( \bdpunit \right) = \left(\begin{array}{ccccc}
    \point & \leftarrow & \point & \rightarrow & G \\
    \point & \mapsfrom & \point & \mapsto & 1
\end{array}\right) \quad \text{and} \quad \mathcal{F}\left( \bdpcounit \right) = \left(\begin{array}{ccccc}
    G & \leftarrow & \point & \rightarrow & \point \\
    1 & \mapsfrom & \point & \mapsto & \point
\end{array}\right) , \]
and hence $Z \left( \bdpunit \right) (1) = [ \{ 1 \} \to G ]$ and $Z \left( \bdpcounit \right) ([X \xrightarrow{f} G]) = [f^{-1}(1)]$.

Now, using \eqref{eq:G_repr_var_when_X_connected} we obtain
\[ [ R_G(N_r) ] [G]^r = [ R_G(N_r, \{ \text{$r + 1$ basepoints} \}) ] = Z \left( \bdpcounit \right) \circ Z \left( \bdpprojplane \right)^r \circ Z \left( \bdpunit \right) (1) . \]
\begin{remark}
    \label{rem:KVar_k_localization}
    To obtain the virtual class $[R_G(N_r)]$ from the equation above, one needs to localize the Grothendieck ring $\K(\Var_k)$ by $[G]$. Consequently, the resulting virtual class will only be known up to an annihilator of $[G]$, but usually this causes no problems when extracting algebraic data.
    For example, the $E$-polynomial of a complex variety can be obtained from its virtual class via the ring morphism $e : \K(\Var_\CC) \to \ZZ[u, v]$. Since $\ZZ[u, v]$ is a domain, as long as $e([G]) \ne 0$, one can obtain the $E$-polynomial of a complex variety from its virtual class in the localized Grothendieck ring.
\end{remark}

Finally, for completeness, we record the TQFT for the bordism $\bdpgenus$, which was shown in \cite{Gonzalez2018,GonzalezLogaresMunoz2020,HablicsekVogel2020} to be given by
\[ Z \left( \bdpgenus \right) : \K(\Var/G) \to \K(\Var/G), \quad \left[ \begin{tikzcd} X \arrow{d}{f} \\ G \end{tikzcd} \right] \mapsto \left[ \begin{tikzcd}[column sep=0em] X \times G^3 \arrow{d} & (x, A, B, C) \arrow[mapsto]{d} \\ G & C f(x) [A, B] C^{-1} \end{tikzcd} \right] . \]

\section{Practically computing the TQFT}
\label{sec:practical}

Even though the construction of the TQFT in the previous section is quite elegant, computing the maps $Z \left( \bdpprojplane \right)$ and $Z \left( \bdpgenus \right)$ is quite cumbersome, and in practice it is more convenient to consider slightly altered maps. These altered maps will have two advantages: they simplify the computations in various aspects, and prevent the need for localization as described in Remark \ref{rem:KVar_k_localization}.

Define the following two morphisms of $\K(\Var_k)$-modules:
\begin{align*}
    Z' \left( \bdprojplane \right) &: \K(\Var/G) \to \K(\Var/G), \quad \left[ \begin{tikzcd} X \arrow{d}{f} \\ G \end{tikzcd} \right] \mapsto \left[ \begin{tikzcd}[column sep=0em, ampersand replacement=\&] X \times G \arrow{d} \& (x, A) \arrow[mapsto]{d} \\ G \& f(x) A^2 \end{tikzcd} \right] \\
    Z' \left( \bdgenus \right) &: \K(\Var/G) \to \K(\Var/G), \quad \left[ \begin{tikzcd} X \arrow{d}{f} \\ G \end{tikzcd} \right] \mapsto \left[ \begin{tikzcd}[column sep=0em, ampersand replacement=\&] X \times G^2 \arrow{d} \& (x, A, B) \arrow[mapsto]{d} \\ G \& f(x) [A, B] \end{tikzcd} \right]
\end{align*}
From the presentations \eqref{eq:representation_variety_orientable_surface} and \eqref{eq:representation_variety_nonorientable_surface} it is clear that
\begin{equation}
    \label{eq:virtual_class_representation_variety_nonorientable_surface}
    [R_G(N_r)] = Z \left( \bdpcounit \right) \circ Z' \left( \bdprojplane \right)^r \circ Z \left( \bdpunit \right)(1)
\end{equation}
and
\begin{equation}
    \label{eq:virtual_class_representation_variety_orientable_surface}
    [R_G(\Sigma_g)] = Z \left( \bdpcounit \right) \circ Z' \left( \bdgenus \right)^g \circ Z \left( \bdpunit \right)(1)
\end{equation}
for any $r \ge 0$ and $g \ge 0$. Note that it is no longer necessary to localize by $[G]$ using these expressions.

Furthermore, we will define the following morphism of $\K(\Var_k)$-modules:
\[ Z' \left( \smallbdmultiplication \right) : \K(\Var/G) \otimes \K(\Var/G) \to \K(\Var/G), \quad \left[ \begin{tikzcd} X \arrow{d}{f} \\ G \end{tikzcd} \right] \otimes \left[ \begin{tikzcd} Y \arrow{d}{g} \\ G \end{tikzcd} \right] \mapsto \left[ \begin{tikzcd}[column sep=-0.25em] X \times Y \arrow{d} & (x, y) \arrow[mapsto]{d} \\ G & f(x) g(y) \end{tikzcd} \right] \]
By definition of these maps, one can express
\begin{equation}
    \label{eq:Z_N_from_multiplication_map}
    Z' \left( \bdprojplane \right) = Z' \left( \smallbdmultiplication \right) \circ \left( \id \otimes \left(Z' \left( \bdprojplane \right) \circ Z \left( \bdpunit \right) \right)(1) \right)
\end{equation}
and similarly
\begin{equation}
    \label{eq:Z_L_from_multiplication_map}
    Z' \left( \bdgenus \right) = Z' \left( \smallbdmultiplication \right) \circ \left( \id \otimes \left(Z' \left( \bdgenus \right) \circ Z \left( \bdpunit \right) \right)(1) \right) .
\end{equation}

\begin{remark}
    Let us explain how these observations simplify the computations. In practice, one needs to evaluate $Z' \left( \bdprojplane \right)$ and $Z' \left( \bdgenus \right)$ on a set of generators (of a submodule) of $\K(\Var/G)$. Hence, for every generator, one needs to deal with the factors $A^2$ and $[A, B]$, respectively, which can be quite hard due to their quadratic nature.   
    Using \eqref{eq:Z_N_from_multiplication_map} and \eqref{eq:Z_L_from_multiplication_map}, these factors only need to be dealt with once in the case where the input is $Z \left( \bdpunit \right) (1)$, which is the easiest case. The price paid is that the multiplication map $Z' \left( \smallbdmultiplication \right)$ should be evaluated on every pair of generators.
    However, those computations are much easier since this map only deals with a single multiplication operation.
    As an additional benefit, using this approach, the computations for the orientable and non-orientable case can mostly be shared.
\end{remark}

\section[AGL1-representation varieties]{$\AGL_1$-representation varieties}
\label{sec:computations_AGL1}

In this section we will prove Theorem \ref{thm:virtual_class_nonorientable_AGL1}, computing the virtual classes of the $G$-representation varieties of non-orientable surfaces in $\K(\Var_k)$ for $G$ equal to the group
\[ \AGL_1 = \left\{ \begin{pmatrix} x & y \\ 0 & 1 \end{pmatrix} \;\Big\vert\; x \ne 0 \right\} \]
of affine transformations of the line, over any field $k$ of characteristic not equal to $2$.
Even though $\K(\Var/G)$ is not finitely generated as $\K(\Var_k)$-module, we would like to express the map $Z' \left( \bdprojplane \right)$ as a matrix with respect to some generators. From the computations that follow, it turns out that this map restricts to the $\K(\Var_k)$-submodule generated by the following elements of $\K(\Var/G)$:
\begin{equation}
    \label{eq:generators_TQFT_AGL1}
    \begin{aligned}
        I &= \textup{ the class of } \left\{ \smatrix{1 & 0 \\ 0 & 1} \right\} , \\
        J &= \textup{ the class of } \left\{ \smatrix{1 & y \\ 0 & 1} \mid y \ne 0 \right\} , \\
        S &= \textup{ the class of } \left\{ \left(\smatrix{x & y \\ 0 & 1}, z \right) \in G \times \AA^1_k \mid z^2 = x \ne 0, 1 \right\} ,
    \end{aligned}
\end{equation}
with all of these varieties being considered naturally as varieties over $G$.

Following the strategy as described in Section \ref{sec:practical}, we first note that
\[ Z \left( \bdpunit \right) (1) = I, \quad Z \left( \bdpcounit \right) (I) = 1, \quad \textup{ and } \quad Z \left( \bdpcounit \right) (J) = Z \left( \bdpcounit \right) (S) = 0 . \]
Next, let us start by computing the virtual class of the morphism $G \to G$ given by $A \mapsto A^2$.
\begin{proposition}
    The virtual class of $G \to G$ given by $A \mapsto A^2$ in $\K(\Var/G)$ is given by
    \[ \left( Z'\left( \bdprojplane \right) \circ Z \left( \bdpunit \right) \right) (1) = (q + 1) I + J + S . \]
\end{proposition}
\begin{proof}
    This morphism can be stratified into the following four strata:
    \begin{align*}
        \left\{ \smatrix{1 & 0 \\ 0 & 1} \right\} \xrightarrow{\quad\sim\quad} I , &\quad \smatrix{1 & 0 \\ 0 & 1} \mapsto \smatrix{1 & 0 \\ 0 & 1} , \\
        \left\{ \smatrix{-1 & y \\ 0 & 1} \right\} \xrightarrow{\quad\sim\quad} \AA^1_k \times I , &\quad \smatrix{-1 & y \\ 0 & 1} \mapsto \left(y, \smatrix{1 & 0 \\ 0 & 1} \right) , \\
        \left\{ \smatrix{1 & y \\ 0 & 1} \mid y \ne 0 \right\} \xrightarrow{\quad\sim\quad} J , &\quad \smatrix{1 & y \\ 0 & 1} \mapsto \smatrix{1 & 2y \\ 0 & 1} , \\
        \left\{ \smatrix{x & y \\ 0 & 1} \mid x \ne 0, \pm 1 \right\} \xrightarrow{\quad\sim\quad} S, &\quad \smatrix{x & y \\ 0 & 1} \mapsto \left(\smatrix{x^2 & (x + 1) y \\ 0 & 1}, x\right) .
    \end{align*}
    Adding up these strata gives the desired result.
\end{proof}

The next step is to compute the image of all pairs of generators in \eqref{eq:generators_TQFT_AGL1} under the map $Z' \left( \smallbdmultiplication \right)$.
\begin{proposition}
    The following equations hold:
    \begin{enumerate}[(i)]
        \item $Z' \left( \smallbdmultiplication \right) (I, X) = X$ for all varieties $X$ over $G$,
        \item $Z' \left( \smallbdmultiplication \right) (J, J) = (q - 1) I + (q - 2) J$,
        \item $Z' \left( \smallbdmultiplication \right) (J, S) = (q - 1) S$,
        \item $Z' \left( \smallbdmultiplication \right) (S, S) = 2 q (q - 3) (I + J) + q (q - 5) S$.
    \end{enumerate}
\end{proposition}
\begin{proof}
    Denote by $\mu : G \times G \to G$ the multiplication map. Clearly $(i)$ holds.
    For $(ii)$, we stratify $J \times J \xrightarrow{\mu} G$ into the two strata
    \begin{align*}
        \left\{ \left( \smatrix{1 & y \\ 0 & 1}, \smatrix{1 & b \\ 0 & 1} \right) \mid b = -y \ne 0 \right\} &\overset{(y, \mu)}{\xrightarrow{\quad\sim\quad}} \left( \AA^1_k \setminus \{ 0 \} \right) \times I , \\
        \left\{ \left( \smatrix{1 & y \\ 0 & 1}, \smatrix{1 & b \\ 0 & 1} \right) \mid 0 \ne b \ne -y \ne 0 \right\} &\overset{(b/y, \mu)}{\xrightarrow{\quad\sim\quad}} \left( \AA^1_k \setminus \{ 0, -1 \} \right) \times J .
    \end{align*}
    For $(iii)$, the map $J \times S \xrightarrow{\mu} G$ is given by
    \[ \left\{ \left( \smatrix{1 & y \\ 0 & 1}, \smatrix{a & b \\ 0 & 1}, c \right) \mid a = c^2 \ne 0, \pm 1 \textup{ and } y \ne 0 \right\} \overset{\left(y, \mu, c\right)}{\xrightarrow{\qquad\sim\qquad}} (\AA^1_k \setminus \{ 0 \}) \times S . \]
    Finally, for $(iv)$ we stratify $S \times S \xrightarrow{\mu} G$ into the three strata
    \begin{align*}
        \left\{ \left( \smatrix{x & y \\ 0 & 1}, z \smatrix{a & b \\ 0 & 1}, c \right) \;\Big\vert\; \substack{x = z^2 \ne 0, 1, \quad b = - y/z^2 \\ \textup{ and } a = c^2 = z^{-2}} \right\} &\overset{\left(cz, y, z, \mu\right)}{\xrightarrow{\qquad\sim\qquad}} \{ \pm 1 \} \times \AA^1_k \times (\AA^1_k \setminus \{ 0, \pm 1 \}) \times I , \\
        \left\{ \left( \smatrix{x & y \\ 0 & 1}, z \smatrix{a & b \\ 0 & 1}, c \right) \;\Big\vert\; \substack{x = z^2 \ne 0, 1, \quad b = 1 - y/z^2 \\ \textup{ and } a = c^2 = z^{-2}} \right\} &\overset{\left(cz, y, z, \mu\right)}{\xrightarrow{\qquad\sim\qquad}} \{ \pm 1 \} \times \AA^1_k \times (\AA^1_k \setminus \{ 0, \pm 1 \}) \times J , \\
        \left\{ \left( \smatrix{x & y \\ 0 & 1}, z \smatrix{a & b \\ 0 & 1}, c \right) \;\Big\vert\; \substack{x = z^2 \ne 0, 1, \quad a = c^2 \ne 0, 1 \\ \textup{ and } c \ne \pm z^{-1}} \right\} &\overset{\left(y, \mu, cz, c \right)}{\xrightarrow{\qquad\sim\qquad}} \AA^1_k \times \left\{ (A, w, c) \in S \times \AA^1_k \mid c \ne 0, \pm 1, \pm w \right\} ,
    \end{align*}
    where the class of the third stratum is seen to be $q (q - 5) S$.
\end{proof}

Now, applying \eqref{eq:Z_N_from_multiplication_map} we obtain
\[ Z' \left( \bdprojplane \right) = \left[\begin{matrix}q + 1 & q - 1 & 2 q \left(q - 3\right)\\1 & 2 q - 1 & 2 q \left(q - 3\right)\\1 & q - 1 & q \left(q - 3\right)\end{matrix}\right] \]
with respect to the generators \eqref{eq:generators_TQFT_AGL1}. Diagonalizing this matrix yields the eigenvalues
\[ 0, \quad q, \quad \textup{ and } \quad q (q - 1) \]
with respective eigenvectors
\[ (q - 3) I + (q - 3) J - S, \quad (q - 1) I - J, \quad \textup{ and } \quad 2 I + 2 J + S . \]
Finally, applying \eqref{eq:virtual_class_representation_variety_nonorientable_surface} proves part $(i)$ of Theorem \ref{thm:virtual_class_nonorientable_AGL1}.

\begin{remark}
    For small values of $r \ge 1$, we have
    \begin{align*}
        [R_{\AGL_1}(N_1)] &= q + 1 , \\
        [R_{\AGL_1}(N_2)] &= 3 q \left(q - 1\right) , \\
        [R_{\AGL_1}(N_3)] &= q^{2} \left(q - 1\right) \left(2 q - 1\right) , \\
        [R_{\AGL_1}(N_4)] &= q^{3} \left(q - 1\right) \left(2 q^{2} - 4 q + 3\right) , \\
        [R_{\AGL_1}(N_5)] &= q^{4} \left(q - 1\right) \left(2 q^{3} - 6 q^{2} + 6 q - 1\right) .
    \end{align*}
\end{remark}

Part $(ii)$ of Theorem \ref{thm:virtual_class_nonorientable_AGL1} can be proven as follows. Any representation in $R_{\AGL_1}(N_r)$ has in the closure of its orbit (under the conjugation action) a unique representation of the form
\[ \left(\begin{pmatrix} \lambda_1 & 0 \\ 0 & 1 \end{pmatrix}, \ldots, \begin{pmatrix} \lambda_r & 0 \\ 0 & 1 \end{pmatrix} \right) \]
with $\lambda_1^2 \cdots \lambda_r^2 = 1$. This shows that the $\AGL_1$-character variety of $N_r$ is given by
\[ X_{\AGL_1}(N_r) \cong \left\{ (\lambda_1, \ldots, \lambda_r) \in \AA^1_k \setminus \{ 0 \} \mid \lambda_1^2 \cdots \lambda_r^2 = 1 \right\} , \]
whose virtual class is $2 (q - 1)^{r - 1}$.

\section[SL2-representation varieties]{$\SL_2$-representation varieties}
\label{sec:computations_SL2}

In this section we will prove Theorem \ref{thm:virtual_class_nonorientable_SL2}, computing the virtual classes of the $G$-representation varieties of non-orientable surfaces in $\K(\Var_k)$ for $G$ equal to the group
\[ \SL_2 = \left\{ \begin{pmatrix} a & b \\ c & d \end{pmatrix} \;\Big\vert\; ad - bc = 1 \right\} \]
of $2 \times 2$ matrices with determinant $1$, over $k = \CC$.
As in the previous section, we will express $Z' \left( \bdprojplane \right)$ as a matrix with respect to the generators of some $\K(\Var_k)$-submodule of $\K(\Var/G)$. It turns out to be convenient to consider the submodule generated by the following elements of $\K(\Var/G)$:
\begin{equation}
    \label{eq:generators_TQFT_SL2}
    \begin{aligned}
        I_+ &= \textup{ the class of } \{ 1 \} , \\
        I_- &= \textup{ the class of } \{ -1 \} , \\
        J_+ &= \textup{ the class of } \left\{ A \in G \mid A \textup{ conjugate to } \left(\begin{smallmatrix} 1 & 1 \\ 0 & 1 \end{smallmatrix}\right) \right\} , \\
        J_- &= \textup{ the class of } \left\{ A \in G \mid A \textup{ conjugate to } \left(\begin{smallmatrix} -1 & 1 \\ 0 & -1 \end{smallmatrix}\right) \right\} , \\
        M &= \textup{ the class of } \{ A \in G \mid \tr(A) \ne \pm 2 \} , \\
        X_2 &= \textup{ the class of } \{ (A, \ell) \in M \times \AA^1_k \mid \ell^2 = \tr(A) - 2 \} , \\
        X_{-2} &= \textup{ the class of } \{ (A, \ell) \in M \times \AA^1_k \mid \ell^2 = \tr(A) + 2 \} , \\
        X_{2, -2} &= \textup{ the class of } \{ (A, \ell) \in M \times \AA^1_k \mid \ell^2 = \tr(A)^2 - 4 \} , \\
        Y &= \textup{ the class of } \{ (A, \ell) \in M \times \AA^1_k \setminus \{ 0 \} \mid \tr(A) = \ell^2 + \ell^{-2} \} ,
    \end{aligned}
\end{equation}
with all of these varieties being considered naturally as varieties over $G$. A useful alternative presentation of the last five generators is as follows:
\begin{align*}
    M &\cong (\GL_2/D \times \AA^1_k \setminus \{ 0, \pm 1 \}) \sslash S_2 \to G, \quad (P, \lambda) \mapsto P \left( \begin{smallmatrix} \lambda & 0 \\ 0 & \lambda^{-1} \end{smallmatrix} \right) P^{-1} \\
    X_2 &\cong (\GL_2/D \times \AA^1_k \setminus \{ 0, \pm 1, \pm i \}) \sslash S_2 \to G, \quad (P, a) \mapsto P \left( \begin{smallmatrix} -a^2 & 0 \\ 0 & -a^{-2} \end{smallmatrix} \right) P^{-1} \\
    X_{-2} &\cong (\GL_2/D \times \AA^1_k \setminus \{ 0, \pm 1, \pm i \}) \sslash S_2 \to G, \quad (P, a) \mapsto P \left( \begin{smallmatrix} a^2 & 0 \\ 0 & a^{-2} \end{smallmatrix} \right) P^{-1} \\
    X_{2, -2} &\cong \GL_2/D \times \AA^1_k \setminus \{ 0, \pm 1 \} \to G, \quad (P, \lambda) \mapsto P \left( \begin{smallmatrix} \lambda & 0 \\ 0 & \lambda \end{smallmatrix} \right) P^{-1} \\
    Y &\cong \GL_2/D \times \AA^1_k \setminus \{ 0, \pm 1, \pm i \} \to G, \quad (P, \lambda) \mapsto P \left( \begin{smallmatrix} x^2 & 0 \\ 0 & x^{-2} \end{smallmatrix} \right) P^{-1}
\end{align*}
where $S_2$ acts on $\GL_2 / D$ by $P \mapsto P \left( \begin{smallmatrix} 0 & 1 \\ 1 & 0 \end{smallmatrix} \right)$, and acts on the coordinates $\lambda$ and $a$ by $\lambda \mapsto \lambda^{-1}$ and $a \mapsto a^{-1}$.
For a better understanding of these elements in the ring $\K(\Var/G)$, we present the following lemma.
\begin{lemma}
    The following relations hold in $\K(\Var/G)$:
    \[ \begin{gathered}
        X_2^2 = 2 X_2, \quad X_{-2}^2 = 2 X_{-2}, \quad X_{2, -2}^2 = 2 X_{2, -2} \\
        \textup{ and } \quad Y = X_2 X_{-2} = X_2 X_{2, -2} = X_{-2} X_{2, -2} .
    \end{gathered} \]
\end{lemma}
\begin{proof}
    The first equality follows from
    \[ X_2 \times_M X_2 = \{ (A, \ell_1, \ell_2) \in M \times \AA^2_k \mid \ell_1^2 = \tr(A) - 2 \textup{ and } \ell_2 = \pm \ell_1 \} \cong X_2 \sqcup X_2 , \]
    and similarly for the second and third. The final two equalities follow from the fact that, if $\ell_1^2 = \tr(A) - 2$ and $\ell_2^2 = \tr(A) + 2$, then $(\ell_1 \ell_2)^2 = \tr(A)^2 - 4$. Finally, the fourth equality follows from the isomorphism
    \[ Y \xrightarrow{\sim} X_2 \times_M X_{-2} = \{ (A, \ell_1, \ell_2) \in M \times \AA^2_k \mid \ell_1^2 = \tr(A) - 2 \textup{ and } \ell_2^2 = \tr(A) + 2 \} \]
    given by $(A, \ell) \mapsto (A, \ell - \ell^{-1}, \ell + \ell^{-1})$ with inverse $(A, \ell_1, \ell_2) \mapsto (A, \tfrac{1}{2} (\ell_1 + \ell_2))$.
\end{proof}

Following the strategy as described in Section \ref{sec:practical}, let us start by computing the class of the morphism $G \to G$ given by $A \mapsto A^2$.
\begin{lemma}
    The virtual class of $G \to G$ given by $A \mapsto A^2$ in $\K(\Var/G)$ equals
    \[ \left( Z'\left( \bdprojplane \right) \circ Z \left( \bdpunit \right) \right) (1) = 2 I_+ + q (q + 1) I_- + 2 J_+ + X_{-2} . \]
\end{lemma}
\begin{proof}
    Write $A = \left( \begin{smallmatrix} a & b \\ c & d \end{smallmatrix} \right)$ and stratify based on the conjugacy class of $A^2$.
    \begin{itemize}
        \item If $A^2 = 1$, then $A = \pm 1$, so we obtain $2 I_+$.
        \item Suppose $A^2 = -1$. If $b = 0$, then $d = a^{-1}$ with $a = \pm i$, contributing $2q I_-$. If $b \ne 0$ and $c = 0$, then $d = a^{-1}$ with $a \ne \pm i$, contributing $2 (q - 1) I_-$. If $b, c \ne 0$, then $a^2 - bc = -1$, contributing $(q - 2) (q - 1) I_-$. In total, we obtain $q (q + 1) I_-$.
        \item Suppose $A^2 \in J_+$. Conjugate to the case $A^2 = \smatrix{1 & 1 \\ 0 & 1}$. Note that there are no solutions for $c \ne 0$. For $c = 0$, solve for $a = d = b/2 = \pm 1$ in order to obtain $2 J_+$.
        \item There are no solutions with $A^2 \in J_-$.
        \item Suppose $A^2 \in M$. This stratum is given by $\left( \GL_2 / D \times \left\{ A \in G \mid A^2 = \left( \begin{smallmatrix} \lambda & 0 \\ 0 & \lambda^{-1} \end{smallmatrix} \right) \right\} \right) \sslash S_2$. There are only solutions with $b = c = 0$ and $d = a^{-1}$ and $a = \lambda^2$. Since $S_2$ acts on $a$ via $a \mapsto d = a^{-1}$, we obtain $X_{-2}$. \qedhere
    \end{itemize}
\end{proof}

The next step is to compute the image of all pairs of generators in \eqref{eq:generators_TQFT_SL2} under the map $Z' \left( \smallbdmultiplication \right)$. As these computations are quite long and repetitive, these have been written down in Appendix \ref{sec:appendix_SL2}.
\begin{remark}
    Except for Lemma \ref{lemma:the_lemma_using_conics}, all of the computations in the appendix are performed in the Grothendieck ring of varieties $\K(\Var/G)$. However, for this particular lemma we pass to the ring $\K^{\PP^1}(\Var/G)$, which is the quotient of $\K(\Var/G)$ by all relations of the form $[P] = [\PP^1_k] [X]$ with $P \to X$ a $\PP^1$-fibration over $G$. As a result, the virtual class $[R_G(N_r)]$ will only be known in $\K^{\PP^1}(\Var_k)$, but usually this causes no problems when extracting algebraic data. For example, the $E$-polynomial of any $\PP^1$-fibration $P \to X$ of complex varieties is known to satisfy $e(P) = e(\PP^1_\CC) e(X)$.
\end{remark}
Applying \eqref{eq:Z_N_from_multiplication_map}, we obtain
\[ Z' \left( \bdprojplane \right) = \scalebox{0.715}{$\footnotesize \left[\begin{matrix}2 & q^{2} + q & 2 q^{2} - 2 & 0 & q^{3} - 3 q^{2} - 2 q & q^{3} - 4 q^{2} - 5 q & 2 q^{3} - 6 q^{2} - 4 q & q^{3} - 4 q^{2} - 5 q & 2 q^{3} - 8 q^{2} - 10 q\\q^{2} + q & 2 & 0 & 2 q^{2} - 2 & q^{3} - 3 q^{2} - 2 q & 2 q^{3} - 6 q^{2} - 4 q & q^{3} - 4 q^{2} - 5 q & q^{3} - 4 q^{2} - 5 q & 2 q^{3} - 8 q^{2} - 10 q\\2 & 0 & q^{2} - q - 2 & 2 q^{2} & q^{3} - 3 q^{2} & q^{3} - 3 q^{2} - 2 q & q^{3} - 4 q^{2} - q & q^{3} - 3 q^{2} - 2 q & q^{3} - 4 q^{2} - 5 q\\0 & 2 & 2 q^{2} & q^{2} - q - 2 & q^{3} - 3 q^{2} & q^{3} - 4 q^{2} - q & q^{3} - 3 q^{2} - 2 q & q^{3} - 3 q^{2} - 2 q & q^{3} - 4 q^{2} - 5 q\\0 & 0 & q^{2} - 1 & q^{2} - 1 & q^{3} - 2 q^{2} - q + 2 & q^{3} - 3 q^{2} - q + 3 & q^{3} - 3 q^{2} - q + 3 & q^{3} - 3 q^{2} - q + 3 & q^{3} - 5 q^{2} - q + 5\\0 & 1 & 0 & - q - 1 & q & - q^{2} + q & q^{2} + q & 0 & 0\\1 & 0 & - q - 1 & 0 & q & q^{2} + q & - q^{2} + q & 0 & 0\\0 & 0 & 1 - q & 1 - q & 2 q - 2 & - q^{2} + 4 q - 3 & - q^{2} + 4 q - 3 & q^{2} + 2 q - 3 & - q^{2} + 6 q - 5\\0 & 0 & q & q & - 2 q & q^{2} - 4 q & q^{2} - 4 q & - 2 q & 2 q^{2} - 6 q\end{matrix}\right]$} \]
with respect to the generators \eqref{eq:generators_TQFT_SL2}. Diagonalizing this matrix yields the eigenvalues
\[ \begin{gathered} 0, \quad -q (q - 1), \quad q (q - 1), \quad q (q - 1) (q + 1), \quad (q - 1) (q + 1), \\ -q (q + 1), \quad 2q (q + 1), \quad 2q (q - 1), \quad q (q + 1) \end{gathered} \]
with respective eigenvectors
\begin{equation}
    \label{eq:eigenvectors_TQFT_SL2}
    \scalebox{0.9}{$\footnotesize \left[\begin{matrix}0\\0\\0\\0\\2\\-1\\-1\\-1\\1\end{matrix}\right]\quad\left[\begin{matrix}- q - 1\\q + 1\\-1\\1\\0\\0\\0\\0\\0\end{matrix}\right]\quad\left[\begin{matrix}- q^{2} + 4 q + 5\\- q^{2} + 4 q + 5\\5 - q\\5 - q\\0\\0\\0\\0\\2\end{matrix}\right]\quad\left[\begin{matrix}1\\1\\1\\1\\1\\0\\0\\0\\0\end{matrix}\right]\quad\left[\begin{matrix}q\\q\\0\\0\\-1\\0\\0\\1\\0\end{matrix}\right]\quad\left[\begin{matrix}q^{2} - 2 q + 1\\- q^{2} + 2 q - 1\\1 - q\\q - 1\\0\\-2\\2\\0\\0\end{matrix}\right]\quad\left[\begin{matrix}1 - q\\q - 1\\1\\-1\\0\\-1\\1\\0\\0\end{matrix}\right]\quad\left[\begin{matrix}q + 1\\q + 1\\1\\1\\0\\0\\0\\-1\\1\end{matrix}\right]\quad\left[\begin{matrix}- q^{2} + 2 q - 1\\- q^{2} + 2 q - 1\\q - 1\\q - 1\\0\\-2\\-2\\0\\2\end{matrix}\right]$} .
\end{equation}
Finally, applying \eqref{eq:virtual_class_representation_variety_nonorientable_surface} proves part $(i)$ of Theorem \ref{thm:virtual_class_nonorientable_SL2}.

\begin{remark}
    The first eigenvector, with eigenvalue $0$, corresponds to the element $2M + Y - X_{2} - X_{-2} - X_{2, -2} \in \K(\Var/G)$. However, this element is not equal to zero, that is, $M + M + Y \ne X_{2} + X_{-2} + X_{2, -2}$ in $\K(\Var/M)$, which can be seen after passing to $\K(\Var_K)$, where $K$ is the function field of $M$. Namely, the left-hand side has two $K$-rational points (corresponding to the two copies of $M$) whereas the right-hand side has no $K$-rational points (neither $X_2, X_{-2}$ nor $X_{2, -2}$ has a section on a dense open of $M$). On the other hand, it is not surprising to encounter this element in the kernel of $Z'\left(\bdprojplane\right)$ as the monodromy representation of $M + M + Y$ agrees with that of $X_2 + X_{-2} + X_{2, -2}$: both are equal to $T^{\oplus 3} \oplus S_2 \oplus S_{-2} \oplus S_2 \otimes S_{-2}$, with $T$ the trivial representation of $\pi_1(M)$ and $S_{\pm 2}$ the $1$-dimensional representations which send a loop around $\tr = \pm 2$ to $-1$, respectively.
\end{remark}

\begin{remark}
    For small values of $r \ge 1$, we have
    \begin{align*}
        [R_{\SL_2}(N_1)] &= 2 , \\
        [R_{\SL_2}(N_2)] &= q \left(q - 1\right) \left(q + 1\right) \left(q + 4\right) , \\
        [R_{\SL_2}(N_3)] &= q \left(q - 1\right) \left(q + 1\right) \left(q^{3} - 8 q - 1\right) , \\
        [R_{\SL_2}(N_4)] &= q \left(q - 1\right) \left(q + 1\right) \left(q^{6} + q^{5} + 13 q^{4} + 3 q^{3} + 13 q^{2} + 1\right) , \\
        [R_{\SL_2}(N_5)] &= q \left(q - 1\right) \left(q + 1\right) \left(q^{9} - 3 q^{7} - 90 q^{5} - 6 q^{4} - 32 q^{3} + 3 q^{2} - 1\right) .
    \end{align*}
\end{remark}
Let us turn for a moment to the case of orientable surfaces. The following proposition describes the virtual class of the morphism $G^2 \to G$ given by $(A, B) \mapsto [A, B]$. This computation was done in \cite{MartinezMunoz2016} on the level of $E$-polynomials, and with some slight modifications this computation can be lifted to the level of virtual classes. Since the case of orientable surfaces is not the main point of this paper, we will omit the explicit details.

\begin{proposition}
    The virtual class of $G^2 \to G$ given by $(A, B) \mapsto [A, B]$ in $\K(\Var/G)$ is given by
    \begin{align*}
        \left( Z'\left( \bdgenus \right) \circ Z( \bdpunit ) \right) (1)
            &= q (q - 1) (q + 1) (q + 4) I_+
            + q (q - 1) (q + 1) I_- \\
            & \quad + q (q - 3) (q + 1) J_+
            + q^2 (q + 3) J_-
            + (q - 1)^2 (q + 1) M \\
            & \quad + 2 q (q + 1) X_{2}
            - q (q + 1) X_{-2}
            - (q - 1)^2 X_{2, -2}
            + q (q - 2) Y . \tag*{\qed}
    \end{align*}
\end{proposition}

Applying \eqref{eq:Z_L_from_multiplication_map}, we obtain
\[ \begin{gathered} Z' \left( \bdgenus \right) = \scalebox{0.715}{\footnotesize $\left[\begin{matrix}q \left(q - 1\right) \left(q + 1\right) \left(q + 4\right) & q \left(q - 1\right) \left(q + 1\right) & q \left(q - 3\right) \left(q - 1\right) \left(q + 1\right)^{2} & q^{2} \left(q - 1\right) \left(q + 1\right) \left(q + 3\right) & q \left(q - 1\right) \left(q + 1\right) \left(q^{3} - 2 q^{2} - 3 q - 2\right)\\q \left(q - 1\right) \left(q + 1\right) & q \left(q - 1\right) \left(q + 1\right) \left(q + 4\right) & q^{2} \left(q - 1\right) \left(q + 1\right) \left(q + 3\right) & q \left(q - 3\right) \left(q - 1\right) \left(q + 1\right)^{2} & q \left(q - 1\right) \left(q + 1\right) \left(q^{3} - 2 q^{2} - 3 q - 2\right)\\q \left(q - 3\right) \left(q + 1\right) & q^{2} \left(q + 3\right) & q \left(q + 1\right) \left(q^{3} + 3\right) & q^{2} \left(q^{3} - 3 q - 6\right) & q^{2} \left(q - 1\right) \left(q^{3} - q^{2} - 4 q - 3\right)\\q^{2} \left(q + 3\right) & q \left(q - 3\right) \left(q + 1\right) & q^{2} \left(q^{3} - 3 q - 6\right) & q \left(q + 1\right) \left(q^{3} + 3\right) & q^{2} \left(q - 1\right) \left(q^{3} - q^{2} - 4 q - 3\right)\\\left(q - 1\right)^{2} \left(q + 1\right) & \left(q - 1\right)^{2} \left(q + 1\right) & q \left(q - 1\right)^{2} \left(q + 1\right)^{2} & q \left(q - 1\right)^{2} \left(q + 1\right)^{2} & \left(q - 1\right)^{2} \left(q + 1\right) \left(q^{3} - q^{2} - 2 q - 2\right)\\2 q \left(q + 1\right) & - q \left(q + 1\right) & - 2 q \left(q + 1\right)^{2} & q \left(q + 1\right)^{2} & q^{2} \left(q + 1\right)\\- q \left(q + 1\right) & 2 q \left(q + 1\right) & q \left(q + 1\right)^{2} & - 2 q \left(q + 1\right)^{2} & q^{2} \left(q + 1\right)\\- \left(q - 1\right)^{2} & - \left(q - 1\right)^{2} & - 2 q \left(q - 1\right)^{2} & - 2 q \left(q - 1\right)^{2} & 2 \left(q - 1\right)^{2} \left(2 q + 1\right)\\q \left(q - 2\right) & q \left(q - 2\right) & 2 q \left(q^{2} - q + 1\right) & 2 q \left(q^{2} - q + 1\right) & - 2 q^{2} \left(2 q - 1\right)\end{matrix}\right.$} \\[10pt]
\scalebox{0.715}{\footnotesize $\left.\begin{matrix}q \left(q - 1\right) \left(q + 1\right) \left(q + 3\right) \left(q^{2} - 3 q - 1\right) & q \left(q - 1\right) \left(q + 1\right)^{2} \left(q^{2} - 4 q - 3\right) & q \left(q - 1\right) \left(q + 1\right) \left(q^{3} - 2 q^{2} - 8 q - 3\right) & q \left(q - 5\right) \left(q - 1\right) \left(q + 1\right) \left(q^{2} + 4 q + 1\right)\\q \left(q - 1\right) \left(q + 1\right)^{2} \left(q^{2} - 4 q - 3\right) & q \left(q - 1\right) \left(q + 1\right) \left(q + 3\right) \left(q^{2} - 3 q - 1\right) & q \left(q - 1\right) \left(q + 1\right) \left(q^{3} - 2 q^{2} - 8 q - 3\right) & q \left(q - 5\right) \left(q - 1\right) \left(q + 1\right) \left(q^{2} + 4 q + 1\right)\\q^{2} \left(q - 1\right) \left(q^{3} - 2 q^{2} - 6 q - 9\right) & q^{2} \left(q - 1\right) \left(q^{3} - 2 q^{2} - 3 q - 6\right) & q^{2} \left(q - 1\right) \left(q^{3} - q^{2} - 6 q - 6\right) & q^{2} \left(q - 5\right) \left(q - 1\right) \left(q^{2} + 2 q + 3\right)\\q^{2} \left(q - 1\right) \left(q^{3} - 2 q^{2} - 3 q - 6\right) & q^{2} \left(q - 1\right) \left(q^{3} - 2 q^{2} - 6 q - 9\right) & q^{2} \left(q - 1\right) \left(q^{3} - q^{2} - 6 q - 6\right) & q^{2} \left(q - 5\right) \left(q - 1\right) \left(q^{2} + 2 q + 3\right)\\\left(q - 1\right)^{2} \left(q + 1\right) \left(q^{3} - 2 q^{2} - 4 q - 3\right) & \left(q - 1\right)^{2} \left(q + 1\right) \left(q^{3} - 2 q^{2} - 4 q - 3\right) & \left(q - 3\right) \left(q - 1\right)^{2} \left(q + 1\right)^{3} & \left(q - 5\right) \left(q - 1\right)^{2} \left(q + 1\right)^{3}\\2 q^{3} \left(q + 1\right) & - 2 q^{2} \left(q - 1\right) \left(q + 1\right) & 0 & 0\\- 2 q^{2} \left(q - 1\right) \left(q + 1\right) & 2 q^{3} \left(q + 1\right) & 0 & 0\\- \left(q - 1\right)^{2} \left(2 q^{2} - 7 q - 3\right) & - \left(q - 1\right)^{2} \left(2 q^{2} - 7 q - 3\right) & \left(q - 1\right)^{2} \left(q^{2} + 8 q + 3\right) & - \left(q - 5\right) \left(q - 1\right)^{2} \left(3 q + 1\right)\\q^{2} \left(2 q^{2} - 11 q + 7\right) & q^{2} \left(2 q^{2} - 11 q + 7\right) & - 6 q^{2} \left(q - 1\right) & 4 q^{2} \left(q - 4\right) \left(q - 1\right)\end{matrix}\right]$}.
\end{gathered} \]
Interestingly, the same set of eigenvectors \eqref{eq:eigenvectors_TQFT_SL2} can be used to diagonalize this matrix, with corresponding eigenvalues
\[ \begin{gathered}
    0, \quad q^2 (q - 1)^2, \quad q^2 (q - 1)^2, \quad q^2 (q - 1)^2 (q + 1)^2, \quad (q - 1)^2 (q + 1)^2, \\ q^2 (q + 1)^2, \quad 4 q^2 (q + 1)^2, \quad 4 q^2 (q - 1)^2, \quad q^2 (q + 1)^2 .
\end{gathered} \]
On the other hand, the fact that both matrices can be simultaniously diagonalized is not too surprising considering the fact that
\[ \bdprojplane \; \circ \; \bdgenus = \bdgenus \, \circ \,\; \bdprojplane \]
as bordisms. Furthermore, it can be seen that
\[ Z' \left( \bdprojplane \right)^3 = Z' \left( \bdprojplane \right) \circ Z' \left( \bdgenus \right) \]
which reflects the fact that the connected sum of three projective planes is equal to the connected sum of a projective plane and a torus.
What is remarkable is that the equality
\[ Z' \left( \bdprojplane \right)^2 = Z' \left( \bdgenus \right) \]
holds for $G = \SL_2$ even though it does not for general $G$. For example, it already fails to hold for $G = \GG_m$. Using \eqref{eq:virtual_class_representation_variety_orientable_surface}, this proves the following theorem.
\begin{theorem}
    Let $\Sigma_g$ be a closed surface of genus $g \ge 0$. Then the virtual class of the $\SL_2$-representation variety in $\K^{\PP^1}(\Var_k)$ of $\Sigma_g$ equals
    \[ [R_{\SL_2}(\Sigma_g)] = [R_{\SL_2}(N_{2g})] . \tag*{\qed} \]
\end{theorem}

\begin{remark}
    An explanation of this relation between the orientable and non-orientable case can be given on the level of $E$-polynomials.
    Formulae by Frobenius \cite[(2.3.8)]{HauselRodriguezVillegas2008} \cite[(9)]{FrobeniusSchur1906} express the number of points of the $G$-representation varieties over a finite field $\FF_q$ in terms of the irreducible characters of $G(\FF_q)$:
    \begin{align*}
        |R_G(\Sigma_g)(\FF_q)| &= \frac{1}{|G(\FF_q)|} \sum_{\chi \in \widehat{G(\FF_q)}} \left( \frac{|G(\FF_q)|}{\chi(1)} \right)^{2g - 2} \\
        |R_G(N_r)(\FF_q)| &= \frac{1}{|G(\FF_q)|} \sum_{\chi \in \widehat{G(\FF_q)}} c_\chi^r \left( \frac{|G(\FF_q)|}{\chi(1)} \right)^{r - 2}
    \end{align*}
    where the sums are taken over the irreducible characters of $G(\FF_q)$ and $c_\chi \in \{ -1, 0, 1 \}$ denotes the Schur indicator of $\chi$. 
    For $q \equiv 1 \pmod 4$, any element of $\SL_2(\FF_q)$ is conjugate to its inverse, and hence
    \[ \chi(g) = \chi(g^{-1}) = \overline{\chi(g)} \]
    for all $g \in \SL_2(\FF_q)$ and irreducible characters $\chi$ of $\SL_2(\FF_q)$. This shows that all irreducible characters $\chi$ of $\SL_2(\FF_q)$ are either real or pseudoreal, that is, $c_\chi = \pm 1$. Therefore,
    \[ |R_{\SL_2}(\Sigma_g)(\FF_q)| = |R_{\SL_2}(N_{2g})(\FF_q)| \]
    for all $q \equiv 1 \pmod 4$, and since these quantities are polynomial in $q$, it follows that $e(R_{\SL_2}(\Sigma_g)) = e(R_{\SL_2}(N_{2g}))$ by Katz' theorem \cite[Theorem 6.1.2]{HauselRodriguezVillegas2008}.
\end{remark}



Let us conclude this section by proving part $(ii)$ of Theorem \ref{thm:virtual_class_nonorientable_SL2}. Denote by $R_{\SL_2}^\textup{red}(N_r)$ the closed subvariety of $R_{\SL_2}(N_r)$ consisting of reducible representations, and by $R_{\SL_2}^\textup{irred}(N_r)$ the open subvariety of irreducible representations. Since $R_{\SL_2}^\textup{irred}(N_r)$ contains the closures of all its orbits, we have
\[ [X_{\SL_2}(N_r)] = [R_{\SL_2}^\textup{red}(N_r) \sslash \SL_2] + [R_{\SL_2}^\textup{irred}(N_r) \sslash \SL_2] , \]
so we will compute both terms on the right-hand side.

First, any representation in $R_{\SL_2}^\textup{red}(N_r)$ has in the closure of its orbit a representation of the form
\[ \left( \smatrix{\lambda_1 & 0 \\ 0 & \lambda_1^{-1}}, \ldots, \smatrix{\lambda_r & 0 \\ 0 & \lambda_r^{-1}} \right) \]
with $\lambda_1^2 \cdots \lambda_r^2 = 1$, which is unique up to a simultaneous permutation $\lambda_i \mapsto \lambda_i^{-1}$. Therefore, we find
\[ [R_{\SL_2}^\textup{red}(N_r) \sslash \SL_2] = \left[ \left\{ \left( \smatrix{\lambda_1 & 0 \\ 0 & \lambda_1^{-1}}, \ldots, \smatrix{\lambda_r & 0 \\ 0 & \lambda_r^{-1}} \right) \;\Big\vert\; \lambda_1^2 \cdots \lambda_r^2 = 1 \right\} \sslash S_2 \right] = (q - 1)^{r - 1} + (q + 1)^{r - 1} . \]
Next, the action of $\SL_2$ on $R_{\SL_2}^\textup{irred}(N_r)$ descents to a free and closed action of $\PGL_2$. Hence, we obtain
\[ [R_{\SL_2}^\textup{irred}(N_r) \sslash \SL_2] = \frac{[R_{\SL_2}(N_r)] - [R_{\SL_2}^\textup{red}(N_r)]}{[\PGL_2]} \]
in $\K^{\PP^1}(\Var_k)$, reducing the problem to computing $[R_{\SL_2}^\textup{red}(N_r)]$. In order to compute this virtual class, define the following subvarieties of $R_{\SL_2}(N_r)$:
\begin{align*}
    R_{\SL_2}^\textup{triv}(N_r) &= \Big\{ \left( \pm \smatrix{1 & 0 \\ 0 & 1}, \ldots , \pm \smatrix{1 & 0 \\ 0 & 1} \right) \Big\} \\
    R_{\SL_2}^\textup{diag}(N_r) &= \left\{ \left( \smatrix{\lambda_1 & 0 \\ 0 & \lambda_1^{-1}}, \ldots , \smatrix{\lambda_r & 0 \\ 0 & \lambda_r^{-1}} \right) \mid \lambda_1^2 \cdots \lambda_r^2 = 1 \right\} \\
    R_{\SL_2}^\textup{upp}(N_r) &= \left\{ \left( \smatrix{\lambda_1 & x_1 \\ 0 & \lambda_1^{-1}}, \ldots , \smatrix{\lambda_r & x_r \\ 0 & \lambda_r^{-1}} \right) \in R_{\SL_2}(N_r) \right\}
\end{align*}
whose virtual classes are given by
\[ [R_{\SL_2}^\textup{triv}(N_r)] = 2^r, \quad [R_{\SL_2}^\textup{diag}(N_r)] = 2 (q - 1)^{r - 1} , \]
and, by means of a TQFT very similar to that of Section \ref{sec:computations_AGL1}, one can obtain
\[ [R_{\SL_2}^\textup{upp}(N_r)] = 2^{r - 1} q^{r - 1} ((-1)^{r} + 1) (q - 1) + 2 q^{r - 1} (q - 1)^{r - 1} . \]
Now, the strata of completely reducible representations is given by
\[ R_{\SL_2}^\textup{triv}(N_r) \sqcup \left( \GL_2 / D \times (R_{\SL_2}^\textup{diag}(N_r) \setminus R_{\SL_2}^\textup{triv}(N_r)) \right) \sslash S_2 , \]
whose virtual class is
\[ q (q - 1) (q + 1)^{r - 1} + q (q - 1)^{r - 1} (q + 1) - 2^r (q - 1) (q + 1) , \]
and the strata of reducible, but not completely reducible, representations is given by
\[ \GL_2 / \left\{ \smatrix{x & y \\ 0 & z} \right\} \times \left( R_{\SL_2}^\textup{upp}(N_r) \setminus R_{\SL_2}^\textup{triv}(N_r) \setminus \left( \AA^1_k \times (R_{\SL_2}^\textup{diag}(N_r) \setminus R_{\SL_2}^\textup{triv}(N_r)) \right) \right) , \]
whose virtual class is
\begin{gather*}
    2^{r - 1} q^{r - 1} (q - 1) (q + 1) ((-1)^r + 1)
    + 2^r (q - 1) (q + 1) \\
    - 2 q (q - 1)^{r - 1} (q + 1)
    + 2 q^{r - 1} (q - 1)^{r - 1} (q + 1) .
\end{gather*}
Adding up these two expressions yields
\begin{align*}
    [R_{\SL_2}^\textup{red}(N_r)]
        &= 2^{r - 1} q^{r - 1} (q - 1) (q + 1) ((-1)^r + 1)
        - q (q + 1) (q - 1)^{r - 1} \\
        &\quad + q (q - 1) (q + 1)^{r - 1}
        + 2 q^{r - 1} (q - 1)^{r - 1} (q + 1) .
\end{align*}
Combining everything, part $(ii)$ of Theorem \ref{thm:virtual_class_nonorientable_SL2} follows.

\begin{remark}
    For small values of $r \ge 1$, we have
    \begin{align*}
        [X_{\SL_2}(N_1)] &= 2, \\
        [X_{\SL_2}(N_2)] &= 3 q - 2, \\
        [X_{\SL_2}(N_3)] &= q^{3} - 6 q - 1, \\
        [X_{\SL_2}(N_4)] &= q^{6} + q^{5} + 11 q^{4} + 9 q^{3} - 5 q^{2} + 2 q + 1, \\
        [X_{\SL_2}(N_5)] &= q^{9} - 3 q^{7} - 2 q^{6} - 84 q^{5} - 10 q^{4} - 30 q^{3} + 9 q^{2} - 1 .
    \end{align*}
\end{remark}


\appendix
\renewcommand{\thesection}{A}
\setcounter{theorem}{0}
\section[Multiplication in SL2]{Multiplication in $\SL_2$}
\label{sec:appendix_SL2}

In the following lemmas, we compute the image under the map $Z' \left( \smallbdmultiplication \right)$ of all pairs of generators in \eqref{eq:generators_TQFT_SL2}. Note that we omit certain cases, but those can directly be obtained from the cases we compute. For example, the case $X = I_-$ is rather trivial, and the cases $X = J_-$ or $X = X_{-2}$ can be derived from the cases $X = J_+$ and $X_2$.

Let us first introduce some notation. When computing $Z' \left( \smallbdmultiplication \right)(X, Y)$ for varieties $X$ and $Y$ over $G$, we write $A$ for an element of $X$ and $B = \smatrix{x & y \\ z & w}$ for an element of $Y$.
When $Y$ is of the form $(\GL_2 / D \times \Lambda) \sslash S_2$  for some $\Lambda$ over $\AA^1_k \setminus \{ 0, \pm 1 \}$, we also write $B = P \smatrix{\mu & 0 \\ 0 & \mu^{-1}} P^{-1}$ with $P = \smatrix{\alpha & \beta \\ \gamma & \delta} \in \GL_2 / D$ and $\mu \in \AA^1_k \setminus \{ 0, \pm 1 \}$. Recall that $S_2$ acts on $(P, \mu)$ via $(P, \mu) \mapsto \left( P \smatrix{0 & 1 \\ 0 & 1}, \mu^{-1} \right)$. More specifically, when $\Lambda = \AA^1_k \setminus \{ 0, \pm 1, \pm i \}$, we write $\mu = a^2$ with $a \in \Lambda$. Similarly, when $X$ is of the form $(\GL_2 / D \times \Lambda) \sslash S_2$ for such $\Lambda$, we write $A = Q \smatrix{\rho & 0 \\ 0 & \rho^{-1}} Q^{-1}$ with $Q \in \GL_2 / D$ and $\rho \in \AA^1_k \setminus \{ 0, \pm 1 \}$, and write $\rho = b^2$ when $\Lambda = \AA^1_k \setminus \{ 0, \pm 1, \pm i \}$.

When considering the strata with $AB \in M$, we commonly base change along the double cover $(\GL_2 / D \times \AA^1_k \setminus \{ 0, \pm 1 \}) \to M$, and write $\lambda$ for the coordinate on $\AA^1_k \setminus \{ 0, \pm 1 \}$. The group $S_2$ acts on this double cover by $(P, \lambda) \mapsto \left(P \smatrix{0 & 1 \\ 1 & 0}, \lambda^{-1} \right)$. The reason for this base change is that it allows to pass to the case $AB = \smatrix{\lambda & 0 \\ 0 & \lambda^{-1}}$ using the symmetry of conjugation.

Often it happens that the set of solutions is invariant under conjugation with respect to some subgroup of $\SL_2$, and one can pick representatives of the orbits of this symmetry. When there is such a symmetry with respect to conjugating with the family $\left\{ \smatrix{1 & x \\ 0 & 1} \right\}$, we speak of \emph{$\GG_a$-invariance}, and when there is symmetry with respect to conjugating with the family $\left\{ \smatrix{x & 0 \\ 0 & x^{-1}} \right\}$, we speak of \emph{$\GG_m$-invariance}.

Finally, to avoid confusion between the different $S_2$-actions, we write $S_2^\lambda$, $S_2^\mu$ and $S_2^\rho$ to differentiate between them.

\begin{lemma}
    \begin{align*}
        Z'\left( \smallbdmultiplication \right)(J_+, J_+) = (q + 1)(q - 1) I_+ + (q - 2) J_+ + q J_- + (q + 1) M - X_{2, -2}
    \end{align*}
\end{lemma}
\begin{proof}
    \begin{itemize}
        \item Case $AB = 1$. This gives $A = B^{-1}$ so we obtain $[J_+] I_+ = (q + 1)(q - 1) I_+$.
        \item Case $AB = -1$. There are no solutions as $\tr(A) = 2 \ne -2 = \tr(-B^{-1})$.
        \item Case $AB \in J_+$. Conjugate to $AB = \smatrix{1 & 1 \\ 0 & 1}$ and solve for $A = \smatrix{w - z & x - y \\ -z & x}$. From $\tr(A) = \tr(B) = 2$ follows that $z = 0$ and $x = w = 1$. Furthermore, $y \ne 0, 1$ as $A, B \ne 1$ so we obtain $(q - 2) J_+$.
        \item Case $AB \in J_-$. Conjugate to $AB = \smatrix{-1 & 1 \\ 0 & -1}$ and solve for $A = \smatrix{-w - z & x + y \\ z & -x}$. Solve for $z = -4$, fix $x = 0$ using $\GG_a$-invariance, and solve for $w = 2$ and $y = 1/4$. We obtain $q J_-$.
        \item Case $AB \in M$. Base change and conjugate to $AB = \smatrix{\lambda & 0 \\ 0 & \lambda^{-1}}$. Solve for $A = \smatrix{\lambda w & - \lambda y \\ - z \lambda^{-1} & x \lambda^{-1}}$ with $w = 2 - x$ and $y = (xw - 1) / z$, using that $z \ne 0$ as $\lambda \ne 1$. Since $\tr(A) = 2$ it follows that $x = \frac{2 \lambda}{\lambda + 1}$. Note that $z \overset{S_2^\lambda}{\mapsto} y = \frac{-(\lambda - 1)^2}{z (\lambda + 1)^2}$, so use $\GG_m$-invariance to fix $z = \frac{i(\lambda - 1)}{(\lambda + 1)}$. Computing the quotient gives $(q + 1) M - X_{2, -2}$. \qedhere
    \end{itemize}
\end{proof}

\begin{lemma}
    \begin{align*}
        Z'\left( \smallbdmultiplication \right)(J_+, M) = q(q - 2) (J_+ + J_-) + (q - 3)(q + 1) M + 2 X_{2, -2}
    \end{align*}
\end{lemma}
\begin{proof}
    Note that $Z'(\smallbdmultiplication)(X, G) = [X] \cdot G$ for all $X \in \K(\Var/G)$. Since $G = I_+ + I_- + J_+ + J_- + M$, the equality can be derived from the earlier lemmas.
\end{proof}

\begin{lemma}
    \begin{align*}
        Z'\left( \smallbdmultiplication \right)(J_+, X_2) = q(q - 3)(J_+ + J_-) + (q - 3)(q + 1)M - (q + 1) X_2 - (q - 3) X_{2, -2} + q Y
    \end{align*}
\end{lemma}
\begin{proof}
    \begin{itemize}
        \item Case $AB = \pm 1$. No solutions.
        \item Case $AB \in J_+$. Conjugate to $AB = \smatrix{1 & 1 \\ 0 & 1}$ and solve for $A = \smatrix{w - z & x - y \\ -z & x}$. Solve for $z = x + w + 2$ and $y = (xw - 1) / z$ and $w = \ell^2 - x + 2$, with $\ell \ne 0, \pm 2i$. We obtain $q (q - 3) J_+$.
        \item Case $AB \in J_-$. Similarly we obtain $q (q - 3) J_-$.
        \item Case $AB \in M$. Base change and conjugate to $AB = \smatrix{\lambda & 0 \\ 0 & \lambda^{-1}}$.
        \begin{itemize}
            \item Case $y = z = 0$. No solutions.
            \item Case $y = 0$ or $z = 0$, but not both. Choose the $S_2^\lambda$-orbit with $z = 0$, solve for $x = w^{-1} = \lambda$ and fix $y = 1$ using $\GG_m$-invariance. We obtain $(q - 1) Y$.
            \item Case $yz \ne 0$. Fix $z = 1$ using $\GG_m$-invariance, solve for $y = xw - 1$ and $w = \lambda^{-1}(2 - x\lambda^{-1})$. Since $y \ne 0$ and $\tr(B) \ne \pm 2$, we have $x \ne \lambda, \frac{2 \lambda}{\lambda + 1}, \frac{-2 \lambda}{\lambda - 1}$. Furthermore, $\ell^2 = \lambda^{-2} (\lambda - 1) (\lambda x - 2 \lambda + x)$. We obtain $(q - 3)(q + 1) M - (q + 1) X_2 - (q - 3) X_{2, -2} + Y$. \qedhere
        \end{itemize}
    \end{itemize}
\end{proof}

\begin{lemma}
    \begin{align*}
        Z'\left( \smallbdmultiplication \right)(J_+, X_{2, -2}) = q (q - 3) (J_+ + J_-) + (q - 3)(q + 1) M + 2 X_{2, -2}
    \end{align*}
\end{lemma}
\begin{proof}
    \begin{itemize}
        \item Case $AB = \pm 1$. No solutions.
        \item Case $AB \in J_+$. Conjugate to $AB = \smatrix{1 & 1 \\ 0 & 1}$. We can choose $\gamma = 1$, $\alpha = 0$ and $\beta = 1$ by $\GG_a$-invariance and lifting $P$ to $\GL_2$. Solve for $\delta = \frac{-(\mu - 1)}{\mu + 1}$. So we obtain $q (q - 3) J_+$.
        \item Case $AB \in J_+$. Similarly we obtain $q (q - 3) J_-$.
        \item Case $AB \in M$. Base change and conjugate to $AB = \smatrix{\lambda & 0 \\ 0 & \lambda^{-1}}$.
        \begin{itemize}
            \item Case $\alpha \gamma = 0$. Choose the $S_2^\lambda$-orbit with $\gamma = 0$ and fix $\alpha = \delta = 1$ by lifting $P$ to $\GL_2$. Then $\tr(A) = 2$ implies $\mu = \lambda$. Furthermore, $\beta$ can be anything, so we obtain $q X_{2, -2}$.
            \item Case $\alpha \gamma \ne 0$. Fix $\alpha = \gamma = 1$ by lifting $P$ to $\GL_2$. Substitute $\beta = u + v$ and $\delta = u - v$. Then $v \ne 0$ and $v \overset{S_2^\lambda}{\mapsto} -v$, and we solve for $u$ using $\tr(A) = 2$. Computing the quotient gives $(q - 3) (q + 1) M - (q - 3) X_{2, -2}$. However, we must subtract the case where $A = 1$, i.e.\ a term $X_{2, -2}$. So we obtain $(q - 3)(q + 1) M - (q - 2) X_{2, -2}$. \qedhere
        \end{itemize}
    \end{itemize}
\end{proof}

\begin{lemma}
    \begin{align*}
        Z'\left( \smallbdmultiplication \right)(J_+, Y) = q (q - 5) (J_+ + J_-) + (q - 5)(q + 1) M - (q - 5) X_{2, -2} + (q - 1) Y
    \end{align*}
\end{lemma}
\begin{proof}
    \begin{itemize}
        \item Case $AB = \pm 1$. No solutions.
        \item Case $AB \in J_+$. Similarly to the case $A \in J_+$ and $B \in X_{2, -2}$, we obtain $q (q - 5) J_+$.
        \item Case $AB \in J_-$. Similarly we obtain $q (q - 5) J_-$.
        \item Case $AB \in M$. Base change and conjugate to $AB = \smatrix{\lambda & 0 \\ 0 & \lambda^{-1}}$.
        \begin{itemize}
            \item Case $\alpha \gamma = 0$. Choose the $S_2^\lambda$-orbit with $\gamma = 0$ and fix $\alpha = \delta = 1$ by lifting $P$ to $\GL_2$. Then $\tr(A) = 2$ implies $a^2 = \lambda$. Furthermore, $\beta$ can be anything, so we obtain $q Y$.
            \item Case $\alpha \gamma \ne 0$. Fix $\alpha = \gamma = 1$ by lifting $P$ to $\GL_2$. Substitute $\beta = u + v$ and $\delta = u - v$. Then $v \ne 0$ and $v \overset{S_2^\lambda}{\mapsto} -v$, and we solve for $u$ using $\tr(A) = 2$. Computing the quotient gives $(q - 5) (q + 1) M - (q - 5) X_{2, -2}$. However, we must subtract the case where $A = 1$, i.e.\ a term $Y$, so we obtain $(q - 5)(q + 1) M - (q - 5) X_{2, -2} - Y$. \qedhere
        \end{itemize}
    \end{itemize}
\end{proof}

\begin{lemma}
    \begin{align*}
        Z'\left( \smallbdmultiplication \right)(M, M) &= q(q^2 - 2q - 1) (I_+ + I_-) + q(q - 3)(q - 1) (J_+ + J_-) \\ &\quad + (q^3 - 4q^2 + 3q + 4) M - 4 X_{2, -2} \\
        Z'\left( \smallbdmultiplication \right)(M, X_2) &= q(q^2 - 3q - 2) (I_+ + I_-) + q(q - 4)(q - 1) (J_+ + J_-) \\ &\quad + (q^3 - 5q^2 + 2q + 6) M + q (X_2 + X_{-2}) + 2(q - 3) X_{2, -2} - 2q Y \\
        Z'\left( \smallbdmultiplication \right)(M, X_{2, -2}) &= q(q - 3)(q + 1) (I_+ + I_-) + q(q - 3)(q - 1) (J_+ + J_-) \\ &\quad + (q - 3)(q - 2)(q + 1) M - 6 X_{2, -2} \\
        Z'\left( \smallbdmultiplication \right)(M, Y) &= q(q - 5)(q + 1) (I_+ + I_-) + q(q - 5)(q - 1) (J_+ + J_-) \\ &\quad + (q - 5)(q - 2)(q + 1) M + 2(q - 5) X_{2, -2} - 2q Y
    \end{align*}
\end{lemma}
\begin{proof}
    Note that $Z'(\smallbdmultiplication)(G, X) = [X] \cdot G$ for all $X \in \K(\Var/G)$. Since $G = I_+ + I_- + J_+ + J_- + M$, the equalities follow from the earlier lemmas.
\end{proof}

\begin{lemma}
    \begin{align*}
        Z'\left( \smallbdmultiplication \right)(X_{2, -2}, X_{2, -2}) &= 2 q (q - 3) (q + 1) (I_+ + I_-) + q (q - 3) (q - 1) (J_+ + J_-) \\ &\quad + (q - 3)^2 (q + 1) M + (q^2 - 4q - 9) X_{2, -2}
    \end{align*}
\end{lemma}
\begin{proof}
    \begin{itemize}
        \item Case $AB = 1$. Solve for $A = B^{-1}$ to obtain $[X_{2, -2} \times_M X_{2, -2}] I_+ = 2 q (q - 3) (q + 1) I_+$.
        \item Case $AB = -1$. Solve for $A = -B^{-1}$ to obtain $[X_{2, -2} \times_M X_{2, -2}] I_- = 2 q (q - 3) (q + 1) I_-$.
        \item Case $AB \in J_+$. Conjugate to $AB = \smatrix{1 & 1 \\ 0 & 1}$.
        \begin{itemize}
            \item Case $\gamma = 0$. Fix $\alpha = \delta = 1$ by lifting $P$ to $\GL_2$ and fix $\beta = 0$ using $\GG_a$-invariance. Then $\rho = \mu^{\pm 1}$, and in both cases $Q$ is well-defined up to a stabilizer of $\smatrix{\rho & 0 \\ 0 & \rho^{-1}}$, so we obtain $2 q (q - 3) J_+$.
            \item Case $\gamma \ne 0$. Fix $\gamma = 1$, $\alpha = 0$ and $\beta = 1$ by lifting $P$ to $\GL_2$ and $\GG_a$-invariance. Then solve for $\delta = \frac{-\beta (\mu - \rho) (\mu \rho - 1)}{\rho (\mu - 1) (\mu + 1)}$. Hence, we obtain $q (q - 3)^2 J_+$.
        \end{itemize}
        \item Case $AB \in J_-$. Similarly we obtain $q (q - 3) (q - 1) J_-$.
        \item Case $AB \in M$. Base change and conjugate to $AB = \smatrix{\lambda & 0 \\ 0 & \lambda^{-1}}$.
        \begin{itemize}
            \item Case $\alpha \gamma = 0$. Choose the $S_2^\lambda$-orbit with $\gamma = 0$. Fix $\alpha = \delta = 1$ by lifting $P$ to $\GL_2$. Then $\rho = (\lambda \mu^{-1})^{\pm 1}$, and $Q$ is well-defined up to a stabilizer of $\smatrix{\rho & 0 \\ 0 & \rho^{-1}}$. We obtain $2 q (q - 5) X_{2, -2}$.
            \item Case $\alpha \gamma \ne 0$ and $\beta \delta = 0$. Choose the $S_2^\lambda$-orbit with $\beta = 0$. Fix $\alpha = \delta = 1$ by lifting $P$ to $\GL_2$. Then $\rho = (\lambda \mu^{-1})^{\pm 1}$ and $Q$ is well-defined up to a stabilizer of $\smatrix{\rho & 0 \\ 0 & \rho^{-1}}$. We obtain $2 (q - 1) (q - 5) X_{2, -2}$.
            \item Case $\alpha \beta \gamma \delta \ne 0$. Fix $\alpha = \gamma = 1$ by lifting $P$ to $\GL_2$. Solve for $\delta = \frac{\beta (\lambda \mu - \rho) (\lambda \mu \rho - 1)}{(\lambda - \mu \rho) (\lambda \rho - \mu)}$. Note that $\beta \ne \delta$ is automatically satisfied as $\rho \ne \lambda^{\pm 1} \mu^{\pm 1}$. Indeed, there are no solutions for $\rho = \lambda^{\pm 1} \mu^{\pm 1}$. We obtain
            \begin{gather*}
                \left(\{ \mu \ne 0, \pm 1, \; \rho \ne 0, \pm 1, \lambda^{\pm 1} \mu^{\pm 1} \} \times \PGL_2 \right) \sslash S_2^\lambda \\
                = (q - 3)^2 (q + 1) M - (3q^2 - 18q + 19) X_{2, -2} . \qedhere
            \end{gather*}
        \end{itemize}
    \end{itemize}
\end{proof}

\begin{lemma}
    \begin{align*}
        Z'\left( \smallbdmultiplication \right)(X_{2, -2}, Y) &= 2 q (q - 5) (q + 1) (I_+ + I_-) + q (q - 5) (q - 1) (J_+ + J_-) \\ &\quad + (q - 5) (q - 3) (q + 1) M + (q - 5) (q + 3) X_{2, -2} - 4 q Y
    \end{align*}
\end{lemma}
\begin{proof}
    \begin{itemize}
        \item Case $AB = 1$. Solve for $A = B^{-1}$ to obtain $[X_{2, -2} \times_M Y] I_+ = 2 q (q - 5) (q + 1) I_+$.
        \item Case $AB = -1$. Solve for $A = -B^{-1}$ to obtain $[X_{2, -2} \times_M Y] I_- = 2 q (q - 5) (q + 1) I_-$.
        \item Case $AB \in J_+$. Conjugate to $AB = \smatrix{1 & 1 \\ 0 & 1}$.
        \begin{itemize}
            \item Case $\gamma = 0$. Fix $\alpha = \delta = 1$ by lifting $P$ to $\GL_2$ and fix $\beta = 0$ using $\GG_a$-invariance. Then $\rho = \mu^{\pm 1}$, and in both cases $Q$ is well-defined up to a stabilizer of $\smatrix{\rho & 0 \\ 0 & \rho^{-1}}$. We obtain $2 q (q - 5) J_+$.
            \item Case $\gamma \ne 0$. Fix $\gamma = 1$, $\alpha = 0$ and $\beta = 1$ by lifting $P$ to $\GL_2$ and $\GG_a$-invariance. Then solve for $\delta = \frac{-\beta (\mu - \rho) (\mu \rho - 1)}{\rho (\mu - 1) (\mu + 1)}$. Hence, we obtain $q (q - 5) (q - 3) J_+$.
        \end{itemize}
        \item Case $AB \in J_-$. Similarly, we obtain $q (q - 5) (q - 1) J_-$.
        \item Case $AB \in M$.
        \begin{itemize}
            \item Case $\alpha \gamma = 0$. Choose the $S_2^\lambda$-orbit with $\gamma = 0$. Fix $\alpha = \delta = 1$ by lifting $P$ to $\GL_2$. Then $\rho = (\lambda \mu^{-1})^{\pm 1}$, and $Q$ is well-defined up to a stabilizer of $\smatrix{\rho & 0 \\ 0 & \rho^{-1}}$. Since $a^2 \ne 0, \pm 1, \pm \lambda$, we obtain $2 q (q - 5) X_{2, -2} - 4 q Y$.
            \item Case $\alpha \gamma \ne 0$ and $\beta \delta = 0$. Choose the $S_2^\lambda$-orbit with $\beta = 0$. Fix $\alpha = \delta = 1$ by lifting $P$ to $\GL_2$. Then $\rho = (\lambda \mu^{-1})^{\pm 1}$ and $Q$ is well-defined up to a stabilizer of $\smatrix{\rho & 0 \\ 0 & \rho^{-1}}$. Since $a^2 \ne 0, \pm 1, \pm \lambda$, we obtain $2 (q - 1) (q - 5) X_{2, -2} - 4 (q - 1) Y$.
            \item Case $\alpha \beta \gamma \delta \ne 0$. Fix $\alpha = \gamma = 1$ by lifting $P$ to $\GL_2$. Solve for $\delta = \frac{\beta (\lambda \mu - \rho) (\lambda \mu \rho - 1)}{(\lambda - \mu \rho) (\lambda \rho - \mu)}$. Note that $\beta \ne \delta$ is automatically satisfied as $\rho \ne \lambda^{\pm 1} \mu^{\pm 1}$. Indeed, there are no solutions for $\rho = \lambda^{\pm 1} \mu^{\pm 1}$. We obtain
            \begin{gather*}
                \left(\{ a \ne 0, \pm 1, \pm i, \; \rho \ne 0, \pm 1, \lambda^{\pm 1} \mu^{\pm 1} \} \times \PGL_2 \right) \sslash S_2^\lambda \\
                = (q - 5) (q - 3) (q + 1) M - (3q^2 - 20q + 25) X_{2, -2} + 4 (q - 1) Y . \qedhere
            \end{gather*}
        \end{itemize}
    \end{itemize}
\end{proof}

\begin{lemma}
    \begin{align*}
        Z'\left( \smallbdmultiplication \right)(Y, Y) &= 4 q (q - 5) (q + 1) (I_+ + I_-) + q (q - 5) (q - 1) (J_+ + J_-) \\ &\quad + (q - 5)^2 (q + 1) M - (q - 5)^2 X_{2, -2} + 2 q (q - 9) Y
    \end{align*}
\end{lemma}
\begin{proof}
    \begin{itemize}
        \item Case $AB = 1$. Solve for $A = B^{-1}$ to obtain $[Y \times_M Y] I_+ = 4 q (q - 5) (q + 1) I_+$.
        \item Case $AB = -1$. Solve for $A = -B^{-1}$ to obtain $[Y \times_M Y] I_- = 4 q (q - 5) (q + 1) I_-$.
        \item Case $AB \in J_+$. Conjugate to $AB = \smatrix{1 & 1 \\ 0 & 1}$.
        \begin{itemize}
            \item Case $\gamma = 0$. Fix $\alpha = \delta = 1$ by lifting $P$ to $\GL_2$ and fix $\beta = 0$ using $\GG_a$-invariance. Then $b^2 = \mu^{\pm 1}$, and in both cases $Q$ is well-defined up to a stabilizer of $\smatrix{\rho & 0 \\ 0 & \rho^{-1}}$. We obtain $4 q (q - 5) J_+$.
            \item Case $\gamma \ne 0$. Fix $\gamma = 1$, $\alpha = 0$ and $\beta = 1$ by lifting $P$ to $\GL_2$ and $\GG_a$-invariance. Then solve for $\delta = \frac{-\beta (\mu - \rho) (\mu \rho - 1)}{\rho (\mu - 1) (\mu + 1)}$. Hence, we obtain $q (q - 5)^2 J_+$.
        \end{itemize}
        \item Case $AB \in J_-$. Similarly, we obtain $q (q - 5) (q - 1) J_-$.
        \item Case $AB \in M$.
        \begin{itemize}
            \item Case $\alpha \gamma = 0$. Choose the $S_2^\lambda$-orbit with $\gamma = 0$. Fix $\alpha = \delta = 1$ by lifting $P$ to $\GL_2$. Then $b^2 = (\lambda \mu^{-1})^{\pm 1}$, and $Q$ is well-defined up to a stabilizer of $\smatrix{\rho & 0 \\ 0 & \rho^{-1}}$. Since $a^2 \ne 0, \pm 1, \pm \lambda$, we obtain $2 q (q - 9) Y$.
            \item Case $\alpha \gamma \ne 0$ and $\beta \delta = 0$. Choose the $S_2^\lambda$-orbit with $\beta = 0$. Fix $\alpha = \delta = 1$ by lifting $P$ to $\GL_2$. Then $b^2 = (\lambda \mu^{-1})^{\pm 1}$ and $Q$ is well-defined up to a stabilizer of $\smatrix{\rho & 0 \\ 0 & \rho^{-1}}$. Since $a^2 \ne 0, \pm 1, \pm \lambda$, we obtain $2 (q - 9) (q - 1) Y$.
            \item Case $\alpha \beta \gamma \delta \ne 0$. Fix $\alpha = \gamma = 1$ by lifting $P$ to $\GL_2$. Solve for $\delta = \frac{\beta (\lambda \mu - \rho) (\lambda \mu \rho - 1)}{(\lambda - \mu \rho) (\lambda \rho - \mu)}$. Note that $\beta \ne \delta$ is automatically satisfied as $\rho \ne \lambda^{\pm 1} \mu^{\pm 1}$. Indeed, there are no solutions for $\rho = \lambda^{\pm 1} \mu^{\pm 1}$.
            We make the substitution $\tilde{\beta} = \beta \frac{\lambda \mu \rho - 1}{\lambda \rho - \mu}$ so that $\tilde{\beta} \overset{S_2^\lambda}{\mapsto} \tilde{\beta}$ and $\tilde{\beta} \overset{S_2^\mu}{\mapsto} \tilde{\beta}^{-1}$. Now,
            \begin{align*}
                \{ a &\ne 0, \pm 1, \pm i, \; b \ne 0, \pm 1, \pm i, \; b^2 \ne \lambda^{\pm 1} a^{\pm 2}, \; \tilde{\beta} \ne 0 \} \\
                &= \{ a \ne 0, \pm 1, \pm i, \; b \ne 0, \pm 1, \pm i, \; \tilde{\beta} \ne 0 \} - \{ a \ne 0, \pm 1, \pm i, \; b^2 = \lambda^{\pm 1} a^{\pm 2}, \; \tilde{\beta} \ne 0 \} \\
                &= (q - 5)^2 (q - 1) - \{ a \ne 0, \pm 1, \pm i, \; b^2 = a^2 \lambda^{\pm 1}, \; \tilde{\beta} \ne 0 \} + \{ b^2 = \lambda^{\pm 1} a^{\pm 2} = \pm 1, \; \tilde{\beta} \ne 0 \} \\
                &= (q - 5)^2 (q - 1) - 2 (q - 5) (q - 1) S_2^\lambda \{ b^2 = \lambda \} + 8 (q - 1) S_2^\lambda \{ a^2 = \lambda \} \\
                &= (q - 5)^2 (q - 1) - 2 S_2^\lambda (q - 9) (q - 1) \{ a^2 = \lambda \} ,
            \end{align*}
            so combined with $Q$ we obtain $(q - 5)^2 (q + 1) M - (q - 5)^2 X_{2, -2} - 2 (q - 9) (q - 1) Y$. \qedhere
        \end{itemize}
    \end{itemize}
\end{proof}

\begin{remark}
    For the final lemmas, rather than computing $Z' \left( \smallbdmultiplication \right)(X_2, -)$, we will instead compute $Z' \left( \smallbdmultiplication \right)(\tilde{X}_2, -)$, where $\tilde{X}_2$ is the class of
    \begin{align*}
        \left(\PGL_2 \times \AA^1_k \setminus \{ 0, \pm 1, \pm i \} \right) \sslash S_2 \to G, \quad (P, a) \mapsto P \smatrix{-a^2 & 0 \\ 0 & -a^{-2}} P^{-1}
    \end{align*}
    with $S_2$ acting by $a \mapsto a^{-1}$ as usual. This turns out to be slightly more convenient, and the desired result can be obtained from the relation
    \[ \tilde{X}_2 = (q + 1) X_2 - Y , \]
    which can easily be verified, in combination with the already known results on $Z' \left( \smallbdmultiplication \right) (Y, -)$.
\end{remark}


\begin{lemma}
    \begin{align*}
        Z'\left( \smallbdmultiplication \right)(\tilde{X}_2, X_{2, -2}) &= q (q - 5) (q - 1) (q + 1) (I_+ + I_-) + q (q - 1) (q^2 - 4 q + 1) (J_+ + J_-) \\ &\quad + (q - 3) (q - 2) (q - 1) (q + 1) M - 6 (q - 1) X_{2, -2} - 2 q (q - 1) Y
    \end{align*}
\end{lemma}
\begin{proof}
    \begin{itemize}
        \item Case $AB = 1$. Solve for $A = B^{-1}$ to obtain $[\tilde{X}_2 \times_M X_{2, -2}] I_+ = q (q - 5) (q - 1) (q + 1) I_+$.
        \item Case $AB = -1$. Solve for $A = - B^{-1}$ to obtain $[\tilde{X}_2 \times_M X_{2, -2}] I_- = q (q - 5) (q - 1) (q + 1) I_-$.
        \item Case $AB \in J_+$. Conjugate to $AB = \smatrix{1 & 1 \\ 0 & 1}$.
        \begin{itemize}
            \item Case $\gamma \delta = 0$. Choose the $S_2^\mu$-orbit with $\gamma = 0$ and fix $\delta = 1$ by lifting $P$ to $\GL_2$ and fix $\beta = 0$ using $\GG_a$-invariance. Then $\rho = a^{\pm 2}$, and in both cases $Q$ is well-defined up to a stabilizer of $\smatrix{\rho & 0 \\ 0 & \rho^{-1}}$. We obtain $2 q (q - 1) (q - 5) J_+$.
            \item Case $\gamma \delta \ne 0$. Fix $\gamma = 1$ by lifting $P$ to $\GL_2$ and fix $\alpha = 0$ using $\GG_a$-invariance. Solve for $\beta = \frac{-\delta \rho (\mu - 1) (\mu + 1)}{(\mu - \rho) (\mu \rho - 1)}$. This requires $\rho \ne 0, \mu^{\pm 1}$ as $\beta, \delta \ne 0$. Note that $\delta \overset{S_2^\mu}{\mapsto} \delta^{-1}$. As usual, $Q$ is well-defined up to a stabilizer of $\smatrix{\rho & 0 \\ 0 & \rho^{-1}}$. We obtain $q (q - 1) (q^2 - 6 q + 11) J_+$.
        \end{itemize}
        \item Case $AB \in J_-$. Similarly, we obtain $q (q - 1) (q^2 - 4 q + 1) J_-$.
        \item Case $AB \in M$. Base change and conjugate to $\smatrix{\lambda & 0 \\ 0 & \lambda^{-1}}$.
        \begin{itemize}
            \item Case $P$ is (anti-)diagonal. Choose the $S_2^\mu$-orbit with $P$ diagonal. Fix $\alpha = \delta = 1$ by lifting $P$ to $\GL_2$ and $\GG_m$-invariance. Then $\rho = (\lambda \mu^{-1})^{\pm 1}$, so choose the $S_2^\lambda$-orbit with $\rho = \lambda \mu^{-1}$. We obtain $(q - 5) (q - 1) X_{2, -2} - 2 (q - 1) Y$.
            \item Case $P$ has \textit{one} zero. Choose the $S_2^\mu$-orbit with $\alpha \gamma = 0$ and the $S_2^\lambda$-orbit with $\gamma = 0$. Fix $\delta = 1$ by lifting $P$ to $\GL_2$ and $\alpha = 1$ by $\GG_m$-invariance. Then $\rho = (\lambda \mu^{-1})^{\pm 1}$, so we obtain $2 (q - 5) (q - 1)^2 X_{2, -2} - 4 (q - 1)^2 Y$.
            \item Case $P$ has no zeros. Fix $\gamma = 1$ and $\alpha = 1$ by lifting $P$ to $\GL_2$ and $\GG_m$-invariance. Solve for $\delta = \frac{\beta (\lambda \mu - \rho) (\lambda \mu \rho - 1)}{(\lambda - \mu \rho) (\lambda \rho - \mu)}$. Note that $\beta \ne \delta$ is automatically satisfied as $\rho \ne \lambda^{\pm 1} \mu^{\pm 1}$. Indeed, there are no solutions for $\rho = \lambda^{\pm 1} \mu^{\pm 1}$.
            We make the substitution $\tilde{\beta} = \beta \frac{\lambda \mu \rho - 1}{\lambda \rho - \mu}$ so that $\tilde{\beta} \overset{S_2^\lambda}{\mapsto} \tilde{\beta}$ and $\tilde{\beta} \overset{S_2^\mu}{\mapsto} \tilde{\beta}^{-1}$. Now,
            \begin{align*}
                \{ a &\ne 0, \pm 1, \pm i, \; \rho \ne 0, \pm 1, \lambda^{\pm 1} a^{\pm 2}, \; \tilde{\beta} \ne 0 \} \sslash S_2^\mu \\
                &= \{ a \ne 0, \pm 1, \pm i, \; \rho \ne 0, \pm 1, \; \tilde{\beta} \ne 0 \} \sslash S_2^\mu - \left\{ \substack{a \ne 0, \pm 1, \pm i, \\ \rho = \lambda^{\pm 1} a^{\pm 2} \ne 0, \pm 1, \; \tilde{\beta} \ne 0 } \right\} \sslash S_2^\mu \\
                &= (q - 1 - 2 S_2^\mu) (q - 3) (q + 1 - S_2^\mu) \sslash S_2^\mu - \left\{ \substack{a \ne 0, \pm 1, \pm i, \\ \rho = a^2 \lambda^{\pm 1}, \; \tilde{\beta} \ne 0} \right\} + \{ a^2 = \pm \lambda^{\pm 1}, \; \tilde{\beta} \ne 0 \} \\
                &= (q - 3) (q - 2) (q - 1) - (q - 5) (q - 1) S_2^\lambda + 2 (q - 1) S_2^\lambda \{ a^2 = \lambda \} \\
                &= (q - 3) (q - 2) (q - 1) - (q - 5 - 2 \{ a^2 = \lambda \}) (q - 1) S_2^\lambda ,
            \end{align*}
            so combined with $Q$ we obtain $(q - 3) (q - 2) (q - 1) (q + 1) M - (q - 1) (2q^2 - 11 q + 11) X_{2, -2} + 2 (q - 1)^2 Y$. \qedhere
        \end{itemize}
    \end{itemize}
\end{proof}

\begin{lemma}
    \begin{align*}
        Z'\left( \smallbdmultiplication \right)(Y, \tilde{X}_2) &= 2 q (q - 5) (q - 1) (q + 1) (I_+ + I_-) + q^2 (q - 5) (q - 1) (J_+ + J_-) \\
        &\quad + (q - 5) (q - 2) (q - 1) (q + 1) M - (q - 5) (q - 2) (q - 1) X_{2, -2} + q (q - 9) (q - 1) Y
    \end{align*}
\end{lemma}
\begin{proof}
    \begin{itemize}
        \item Case $AB = 1$. Solve for $A = B^{-1}$ to obtain $[Y \times_M \tilde{X}_2] I_+ = 2 q (q - 5) (q - 1) (q + 1) I_+$.
        \item Case $AB = -1$. Solve for $A = -B^{-1}$ to obtain $[Y \times_M \tilde{X}_2] I_- = 2 q (q - 5) (q - 1) (q + 1) I_-$.
        \item Case $AB \in J_+$. Conjugate to $AB = \smatrix{1 & 1 \\ 0 & 1}$.
        \begin{itemize}
            \item Case $\gamma \delta = 0$. Choose the $S_2^\mu$-orbit with $\gamma = 0$ and fix $\delta = 1$ by lifting $P$ to $\GL_2$ and fix $\beta = 0$ using $\GG_a$-invariance. Then $b^2 = a^{\pm 2}$, that is, $b = \pm a^{\pm 1}$. In all cases $Q$ is well-defined up to a stabilizer of $\smatrix{\rho & 0 \\ 0 & \rho^{-1}}$. We obtain $4 q (q - 5) (q - 1) J_+$.
            \item Case $\gamma \delta \ne 0$. Fix $\gamma = 1$ by lifting $P$ to $\GL_2$ and fix $\alpha = 0$ using $\GG_a$-invariance. Solve for $\beta = \frac{-\delta \rho (\mu - 1) (\mu + 1)}{(\mu - \rho) (\mu \rho - 1)}$. This requires $b^2 \ne 0, \pm 1, a^{\pm 2}$ as $\beta, \delta \ne 0$. Note that $\delta \overset{S_2^\mu}{\mapsto} \delta^{-1}$. As usual, $Q$ is well-defined up to a stabilizer of $\smatrix{\rho & 0 \\ 0 & \rho^{-1}}$. We obtain $q (q - 5) (q - 4) (q - 1) J_+$.
        \end{itemize}
        \item Case $AB \in J_-$. Similarly, we obtain $q^2 (q - 5) (q - 1) J_-$.
        \item Case $AB \in M$. Base change and conjugate to $AB = \smatrix{\lambda & 0 \\ 0 & 0 \lambda^{-1}}$.
        \begin{itemize}
            \item Case $P$ is (anti-)diagonal. Choose the $S_2^\mu$-orbit with $P$ diagonal. Fix $\alpha = \delta = 1$ by lifting $P$ to $\GL_2$ and $\GG_m$-invariance. Then $b^2 = (\lambda a^{-2})^{\pm 1}$, so choose the $S_2^\lambda$-orbit with $b^2 = \lambda a^{-2}$. We obtain $(q - 9) (q - 1) Y$.
            \item Case $P$ has \textit{one} zero. Choose the $S_2^\mu$-orbit with $\alpha \gamma = 0$ and the $S_2^\lambda$-orbit with $\gamma = 0$. Fix $\delta = 1$ by lifting $P$ to $\GL_2$ and $\alpha = 1$ by $\GG_m$-invariance. Then $b^2 = (\lambda a^{-2})^{\pm 1}$, so we obtain $2 (q - 9) (q - 1)^2 Y$.
            \item Case $P$ has no zeros. Fix $\gamma = 1$ and $\alpha = 1$ by lifting $P$ to $\GL_2$ and $\GG_m$-invariance. Solve for $\delta = \frac{\beta (\lambda \mu - \rho) (\lambda \mu \rho - 1)}{(\lambda - \mu \rho) (\lambda \rho - \mu)}$. Note that $\beta \ne \delta$ is automatically satisfied as $b^2 \ne \lambda^{\pm 1} a^{\pm 2}$. Indeed, there are no solutions for $b^2 = \lambda^{\pm 1} a^{\pm 2}$.
            We make the substitution $\tilde{\beta} = \beta \frac{\lambda \mu \rho - 1}{\lambda \rho - \mu}$ so that $\tilde{\beta} \overset{S_2^\lambda}{\mapsto} \tilde{\beta}$ and $\tilde{\beta} \overset{S_2^\mu}{\mapsto} \tilde{\beta}^{-1}$. Now,
            \begin{align*}
                \{ a &\ne 0, \pm 1, \pm i, \; b \ne 0, \pm 1, \pm i, \; b^2 \ne \lambda^{\pm 1} a^{\pm 2}, \; \tilde{\beta} \ne 0 \} \sslash S_2^\mu \\
                &= \left\{ \substack{a \ne 0, \pm 1, \pm i, \\ b \ne 0, \pm 1, \pm i, \; \tilde{\beta} \ne 0 } \right\} \sslash S_2^\mu - \left\{ \substack{a \ne 0, \pm 1, \pm i, \\ b^2 = \lambda^{\pm 1} a^{\pm 2} \ne 0, \pm 1, \; \tilde{\beta} \ne 0 } \right\} \sslash S_2^\mu \\
                &= (q - 1 - 2 S_2^\mu) (q - 5) (q + 1 - S_2^\mu) \sslash S_2^\mu - \left\{ \substack{a \ne 0, \pm 1, \pm i, \\ b^2 = \lambda^{\pm 1} a^2, \; \tilde{\beta} \ne 0} \right\} + \left\{ \substack{b^2 = \lambda^{\pm 1} a^2 = \pm 1, \\ \tilde{\beta} \ne 0} \right\} \\
                &= (q - 5) (q - 2) (q - 1) - (q - 5) (q - 1) S_2^\lambda \{ b^2 = \lambda \} + 4 (q - 1) S_2^\lambda \{ a^2 = \lambda \} \\
                &= (q - 5) (q - 2) (q - 1) - (q - 9) (q - 1) S_2^\lambda \{ a^2 = \lambda \} ,
            \end{align*}
            so combined with $Q$ we obtain $(q - 5) (q - 2) (q - 1) (q + 1) M - (q - 5) (q - 2) (q - 1) X_{2, -2} - (q - 9) (q - 1)^2 Y$. \qedhere
        \end{itemize}
    \end{itemize}
\end{proof}

\begin{lemma}
    \label{lemma:the_lemma_using_conics}
    \begin{align*}
        Z'\left( \smallbdmultiplication \right)(\tilde{X}_2, \tilde{X}_2)
            &= 2 q (q - 2) (q - 1)^2 (q + 1) I_+ + q (q - 5) (q - 1)^2 (q + 1) I_- \\
            &\quad + q^3 (q - 5) (q - 1) J_+ + q (q - 1)^2 (q^2 - 3 q - 1) J_- \\
            &\quad + (q - 2)^2 (q - 1)^2 (q + 1) M - q (q - 3) (q + 1)^2 X_{-2} \\
            &\quad - (q - 2)^2 (q - 1)^2 X_{2, -2} + q (q^3 - 6 q^2 + 7 q - 6) Y
    \end{align*}
\end{lemma}
\begin{proof}
    \begin{itemize}
        \item Case $AB = 1$. Solve for $A = B^{-1}$. Then $b^2 = a^{\pm 2}$, that is, $b = \pm a^{\pm 1}$, so choose the $S_2^\rho$-orbit with $b = \pm a^{-1}$. Then $Q = P \smatrix{\xi & 0 \\ 0 & 1}$ with $\xi \ne 0$ and $\xi \overset{S_2^\mu}{\mapsto} \xi^{-1}$. Taking the $S_2^\mu$-quotient, we obtain $2 q (q - 2) (q - 1)^2 (q + 1) I_+$.
        \item Case $AB = -1$. Solve for $A = -B^{-1}$. Then $b^2 = -a^{\pm 2}$, that is, $b = \pm i a^{\pm 1}$, so choose the $S_2^\rho$-orbit with $b = \pm i a^{-1}$, and choose the $S_2^\mu$-orbit with $b = i a^{-1}$. Then $Q$ is well-defined up to a stabilizer of $\smatrix{\rho & 0 \\ 0 & \rho^{-1}}$. We obtain $q (q - 5) (q - 1)^2 (q + 1) I_-$.
        \item Case $AB \in J_+$. Conjugate to $AB = \smatrix{1 & 1 \\ 0 & 1}$.
        \begin{itemize}
            \item Case $\gamma \delta = 0$. Choose the $S_2^\mu$-orbit with $\gamma = 0$. Fix $\delta = 1$ by lifting $P$ to $\GL_2$ and fix $\beta = 0$ using $\GG_a$-invariance. Then $b^2 = a^{\pm 2}$, that is, $b = \pm a^{\pm 1}$, so choose the $S_2^\rho$-orbit with $b = \pm a$. Furthermore, $Q$ is well-defined up to a stabilizer of $\smatrix{\rho & 0 \\ 0 & \rho^{-1}}$. We obtain $2 q (q - 5) (q - 1)^2 J_+$.
            \item Case $\gamma \delta \ne 0$. Fix $\gamma = 1$ by lifting $P$ to $\GL_2$ and fix $\alpha = 0$ using $\GG_a$-invariance. Solve for $\beta = \frac{-\delta \rho (\mu - 1) (\mu + 1)}{(\mu - \rho) (\mu \rho - 1)}$.  This requires $b^2 \ne 0, \pm 1, a^{\pm 2}$ as $\beta, \delta \ne 0$. Note that $\delta \overset{S_2^\mu}{\mapsto} \delta^{-1}$. Solve for $Q = \smatrix{-\xi \rho & -1 \\ \xi (\mu - \rho) & \mu \rho - 1}$ for some $\xi \ne 0$ with $\xi \overset{S_2^\mu}{\mapsto} - \xi$ and $\xi \overset{S_2^\rho}{\mapsto} \xi^{-1}$. We obtain $q (q - 5) (q - 1) (q^2 - 2 q + 2) J_+$.
        \end{itemize}
        \item Case $AB \in J_-$. Conjugate to $AB = \smatrix{-1 & 1 \\ 0 & -1}$.
        \begin{itemize}
            \item Case $\gamma \delta = 0$. Similarly to the above we obtain $2 q (q - 5) (q - 1)^2 J_+$.
            \item Case $\gamma \delta \ne 0$. Fix $\gamma = 1$ by lifting $P$ to $\GL_2$ and fix $\alpha = 0$ using $\GG_a$-invariance. Solve for $\beta = \frac{-\delta \rho (\mu - 1) (\mu + 1)}{(\mu - \rho) (\mu \rho - 1)}$.  This requires $b^2 \ne 0, \pm 1, -a^{\pm 2}$ as $\beta, \delta \ne 0$. Note that the cases $b = \pm i a^{\pm 1}$ behave as the product $S_2^\mu S_2^\rho$ under the action of $S_2^\mu$ and $S_2^\rho$. Furthermore, $\delta \overset{S_2^\mu}{\mapsto} \delta^{-1}$. Solve for $Q = \smatrix{-\xi \rho & -1 \\ \xi (\mu - \rho) & \mu \rho - 1}$ for some $\xi \ne 0$ with $\xi \overset{S_2^\mu}{\mapsto} -\xi$ and $\xi \overset{S_2^\rho}{\mapsto} \xi^{-1}$. We obtain $q (q - 1)^2 (q^2 - 5 q + 9) J_+$.
        \end{itemize}
        \item Case $AB \in M$. Base change and conjugate to $AB = \smatrix{\lambda & 0 \\ 0 & \lambda^{-1}}$.
        \begin{itemize}
            \item Case $P$ is (anti)-diagonal. Choose the $S_2^\mu$-orbit with $P$ diagonal. Fix $\delta = 1$ by lifting $P$ to $\GL_2$ and fix $\alpha = 1$ using $\GG_m$-invariance. Then $b^2 = (\lambda a^{-2})^{\pm 1}$, so choose the $S_2^\rho$-orbit with $b^2 = \lambda a^{-2}$. Furthermore, $Q = \smatrix{\xi & 0 \\ 0 & 1}$ with $\xi \overset{S_2^\lambda}{\mapsto} \xi^{-1}$. We obtain $(q - 3) (q + 1)^2 X_{-2} - (5 q^2 - 12 q + 3) Y$.
            \item Case $P$ has \textit{one} zero. Choose the $S_2^\mu$-orbit with $\alpha \gamma = 0$ and the $S_2^\lambda$-orbit with $\gamma = 0$. Fix $\delta = 1$ by lifting $P$ to $\GL_2$ and fix $\alpha = 1$ using $\GG_m$-invariance. Choose the $S_2^\rho$-orbit with $b^2 = \lambda a^{-2}$. Furthermore, $Q$ is well-defined up to stabilizer of $\smatrix{\rho & 0 \\ 0 & \rho^{-1}}$. We obtain $(q - 9)(q - 1)^3 Y$.
            \item Case $P$ has no zeros. Fix $\gamma = 1$ by lifting $P$ to $\GL_2$ and fix $\alpha = 1$ using $\GG_m$-invariance. Solve for $\delta = \frac{\beta(\lambda \mu - \rho)(\lambda \mu \rho - 1)}{(\lambda - \mu \rho)(\lambda \rho - \mu)}$. Note that $\beta \ne \delta$ is automatically satisfied as $b^2 \ne \lambda^{\pm 1} a^{\pm 2}$. Indeed, there are no solutions for $b^2 = \lambda^{\pm 1} a^{\pm 2}$. The $S_2$-actions on $\beta$ are given by
            \[ \beta \overset{S_2^\lambda}{\mapsto} \frac{\beta(\lambda \mu - \rho)(\lambda \mu \rho - 1)}{(\lambda - \mu \rho)(\lambda \rho - \mu)}, \quad \beta \overset{S_2^\mu}{\mapsto} \beta^{-1}, \quad \beta \overset{S_2^\rho}{\mapsto} \beta . \]
            Solve for $Q = \smatrix{\xi \lambda(\mu - \lambda\rho) & \lambda(\mu\rho - \lambda) \\ \xi (\lambda\mu - \rho) & \lambda\mu\rho - 1}$ for some $\xi \ne 0$. After a substitution $\tilde{\xi} = \tilde{\xi} \frac{\lambda - \mu\rho}{\lambda\rho - \mu}$, the $S_2$-actions on $\tilde{\xi}$ are given by
            \[ \tilde{\xi} \overset{S_2^\lambda}{\mapsto} \tilde{\xi} \frac{(\lambda \rho - \mu)(\lambda \mu \rho - 1)}{(\lambda - \mu\rho)(\lambda \mu - \rho)}, \quad \tilde{\xi} \overset{S_2^\mu}{\mapsto} \tilde{\xi}, \quad \tilde{\xi} \overset{S_2^\rho}{\mapsto} \tilde{\xi}^{-1} . \]
            Explicitly computing the $S_2^\mu$-quotient and the $S_2^\rho$-quotient yields
            \[ \{ x, y \ne 0, \pm 2, \; u, v, c, d \in k \mid v^2 = (x^2 - 4)(u^2 - 4) \text{ and } d^2 = (y^2 - 4)(c^2 - 4) \text{ and } F \ne 0  \} \]
            where $x = a + a^{-1}$, $y = b + b^{-1}$, $u = \beta + \beta^{-1}$, $v = (\beta - \beta^{-1})(a - a^{-1})$, $c = \tilde{\xi} + \tilde{\xi}^{-1}, \eta = (\tilde{\xi} - \tilde{\xi}^{-1})(b - b^{-1})$ and $F = r^2 + s^2 + t^2 - rst - 4$. We make suitable substitutions
            \[ \begin{pmatrix} u \\ v \end{pmatrix} = \begin{pmatrix} \frac{2r - t(x^2 - 2)}{x(x^2 - 4)} & \lambda - \lambda^{-1} \\ \lambda - \lambda^{-1} & \frac{2r - t(x^2 - 2)}{x} \end{pmatrix} \begin{pmatrix} \tilde{u} \\ \tilde{v} \end{pmatrix} \text{ and } \begin{pmatrix} c \\ d \end{pmatrix} = \begin{pmatrix} \frac{2s - t(y^2 - 2)}{y(y^2 - 4} & \lambda - \lambda^{-1} \\ \lambda - \lambda^{-1} & \frac{2s - t(y^2 - 2)}{y} \end{pmatrix} \begin{pmatrix} \tilde{c} \\ \tilde{d} \end{pmatrix} , \]
            so that $\tilde{u}, \tilde{v}$ and $\tilde{c}, \tilde{d}$ are $S_2^\lambda$-invariant, and the relevant equations are (after some rescaling) given by
            \[ \tilde{u}^2 (x^2 - 4) - \tilde{v}^2 = F \quad \text{and} \quad \tilde{c}^2 (y^2 - 4) - \tilde{d}^2 = F . \]
            This is a family of conics over a family of conics over the base $\{ x, y \ne 0, \pm 2, \; F \ne 0 \}$. Note that both families have conics with non-zero discriminant (as $x, y \ne \pm 2$). So we obtain
            \begin{align*}
                (q + 1) &\left( (q + 1) \{ x, y \ne 0, \pm 2, \; F \ne 0 \} - \{ a \ne 0, \pm 1, \pm i, \; y \ne 0, \pm 2, \; F \ne 0 \} \right) \\
                        &- \left( (q + 1) \{ x \ne 0, \pm 2, \; b \ne 0, \pm 1, \pm i, \; F \ne 0 \} - \{ a, b \ne 0, \pm 1, \pm i, \; F \ne 0 \} \right) \\
                        &= (q - 2)^2 (q - 1)^2 - (q - 3)(q + 1)^2 \{ a^2 = \lambda \} + (5q^2 - 12q + 3) S_2^\lambda ,
            \end{align*}
            since
            \begin{align*}
                \{ x, y \ne 0, \pm 2, \; F \ne 0 \} &= (q - 3)^2 - (q - 3 - 3 S_2^\lambda) \{ a^2 = \lambda \} , \\
                \{ x \ne 0, \pm 2, \; b \ne 0, \pm 1, \pm i, \; F \ne 0 \} &= (q - 5)(q - 3) - S_2^\lambda  (q - 9) \{ a^2 = \lambda \} , \\
                \{ a \ne 0, \pm 1, \pm i, \; y \ne 0, \pm 2, \; F \ne 0 \} &= (q - 5)(q - 3) - S_2^\lambda  (q - 9) \{ a^2 = \lambda \} , \\
                \{ a, b \ne 0, \pm 1, \pm i, \; F \ne 0 \} &= (q - 5)^2 - 2 S_2^\lambda (q - 9) \{ a^2 = \lambda \} .
            \end{align*}
            Hence, we obtain $(q - 2)^2 (q - 1)^2 (q + 1) M - (q - 3) (q + 1)^3 X_{-2} - (q - 2)^2 (q - 1)^2 X_{2, -2} + 2 (3 q^3 - 9 q^2 + 5 q - 3) Y$. \qedhere
        \end{itemize}
    \end{itemize}
\end{proof}

\newpage 
\printbibliography

\end{document}